%&amstex          
\input amstex\documentstyle{amsppt}  
\pagewidth{12.5cm}\pageheight{19cm}\magnification\magstep1
\topmatter
\title Irreducible representations of finite spin groups\endtitle
\author G. Lusztig\endauthor
\address{Department of Mathematics, M.I.T., Cambridge, MA 02139}\endaddress
\thanks{Supported in part by the National Science Foundation.}\endthanks
\endtopmatter   
\document
\redefine\O{\Omega}
\define\boa{\bold a}
\define\ufa{\un{\fA}}
\define\tufa{\ti{\ufa}}

\define\tSO{\wt{SO}}
\define\Irr{\text{\rm Irr}}
\define\Wi{\text{\rm Wi}}

\define\uT{\un T}
\define\uA{\un A}
\define\uXi{\un{\Xi}}

\define\tir{\ti{\rho}}

\define\tcs{\ti{\cs}}

\define\dg{\dot g}

\define\uf{\un f}

\define\ha{\hat a}

\define\pe{\perp}
\define\si{\sim}
\define\wt{\widetilde}
\define\sqc{\sqcup}

\define\qua{\quad}

\define\hf{\hat f}

\define\lb{\linebreak}

\define\op{\oplus}

\define\part{\partial}
\define\em{\emptyset}

\define\n{\notin}

\define\m{\mapsto}
\define\do{\dots}

\define\lra{\leftrightarrow}

\define\sub{\subset}    

\define\T{\times}
\define\ti{\tilde}
\define\nl{\newline}
\redefine\i{^{-1}}

\define\un{\underline}
\define\ov{\overline}

\define\bbq{\bar{\QQ}_l}

\define\ad{\text{\rm ad}}
\define\Ad{\text{\rm Ad}}
\define\Hom{\text{\rm Hom}}

\define\card{\text{\rm card}}
\define\bst{\bigstar}
\define\he{\heartsuit}

\define\a{\alpha}
\redefine\b{\beta}

\define\g{\gamma}
\redefine\d{\delta}
\define\e{\epsilon}
\define\et{\eta}
\define\io{\iota}
\redefine\o{\omega}
\define\p{\pi}
\define\ph{\phi}
\define\ps{\psi}
\define\r{\rho}
\define\s{\sigma}
\redefine\t{\tau}
\define\th{\theta}
\define\k{\kappa}
\redefine\l{\lambda}
\define\z{\zeta}
\define\x{\xi}

\define\vp{\varpi}

\redefine\G{\Gamma}
\redefine\D{\Delta}

\define\Ph{\Phi}
\define\Ps{\Psi}
\define\La{\Lambda}

\define\boc{\bold c}

\define\kk{\bold k}

\define\uu{\bold u}

\redefine\AA{\bold A}

\define\FF{\bold F}

\define\NN{\bold N}

\define\QQ{\bold Q}

\define\VV{\bold V}

\define\ZZ{\bold Z}

\define\ca{\Cal A}

\define\cc{\Cal C}

\define\cf{\Cal F}

\define\co{\Cal O}
\define\cp{\Cal P}

\define\cs{\Cal S}
\define\ct{\Cal T}
\define\cu{\Cal U}

\define\cx{\Cal X}

\define\fa{\frak a}
\define\fb{\frak b}
\define\fc{\frak c}
\define\fd{\frak d}
\define\fe{\frak e}

\define\fs{\frak s}

\define\fA{\frak A}

\define\fF{\frak F}

\define\fN{\frak N}

\define\fS{\frak S}

\define\fU{\frak U}

\define\tf{\ti f}
\define\tg{\ti g}

\define\tit{\ti t}

\define\tz{\ti z}

\define\ty{\ti y}
\define\tA{\ti A}

\define\tE{\ti E}

\define\tO{\ti O}
\define\tP{\ti P}

\define\tT{\ti T}

\define\tZ{\ti Z}

\define\bul{\bullet}

\define\tri{\triangle}
\define\DL{DL}
\define\IRC{L1}
\define\WA{L2}
\define\OR{L3}
\define\ICC{L4}
\define\AST{L5}

\head Introduction\endhead
\subhead 0.1\endsubhead
Let $G$ be a connected reductive group defined over a finite field $\FF_q$ and let $F:G@>>>G$
be the corresponding Frobenius map so that the fixed point $G^F$ is a finite group.
In this paper we are interested in the problem of classification of the irreducible 
representations of $G^F$ over $\bbq$, an algebraic closure of  the field of $l$-adic 
numbers ($l$ is a fixed prime number not dividing $q$). When the centre of $G$ is
connected this problem was solved in \cite{\OR} using the $l$-adic cohomology
approach of \cite{\DL}; namely in \cite{\OR} it is shown that in this case the
irreducible representations of $G^F$ can be parametrized in terms of semisimple
conjugacy defined over $\FF_q$ in the dual group of $G$ and by certain 
(essentially combinatorial) objects called unipotent representations and which are
attached to a (usually) smaller group and can themselves be classified.
A generalization of this parametrization to the case where the centre of $G$ is
not necessarily connected was stated without proof in \cite{\WA} and with a partial
proof in \cite{\AST}. Namely in \cite{\AST} it was shown that the statement in
\cite{\WA} can be proved assuming a certain multiplicity one statement for the
finite spin groups in dimension $4k$ in odd characteristic. Moreover it was shown 
how this multiplicity one
statement can be reduced to a computation which however was not given there.
In this paper we perform the needed computation; namely we compute explicitly
(in terms of a generating function) the difference between the number of
conjugacy classes of the two forms of a spin group in dimension $2k$ (see Section 3); 
we compute separately the analogous number coming from counting data in the dual group
(see Section 4). We find
that these two differences coincide. This is sufficient for proving the required
multiplicity one statement (5.4(f)) and provides the missing step for the classification
problem stated above.

\subhead 0.2. Notation\endsubhead
If $f$ is a permutation of a set $S$ we write $S^f=\{x\in S;f(x)=x\}$. If $S$ is finite we
write $|S|$ instead of $\card(S)$. 

If $P$ is a property we set $\d_P=1$ if $P$ is true, $\d_P=0$ if $P$ is false.

If $S$ is a set such that $|S|=2$, there is a unique free transitive action of $\ZZ/2$ on 
$S$ denoted by $n:x\m n+x$ (here $n\in\ZZ/2$, $x\in S$). 

We fix an algebraic closure $\kk$ of a finite field $\FF_q$ with $|\FF_q|=q$, $q$ odd. Let
$\kk^*=\kk-\{0\}$. 

Let $\Wi(\FF_q)$ be the Witt group of $\FF_q$; let $\Wi^0(\FF_q)$ be the subgroup of 
$\Wi(\FF_q)$ generated by classes of quadratic forms of even dimension. We have naturally
$\Wi^0(\FF_q)=\ZZ/2$. 

All algebraic groups are assumed to be affine and over $\kk$. For an algebraic group $G$ we
denote by $G^0$ the identity component of $G$ and we set $\ov{G}=G/G^0$. For $g\in G$ we
denote by $g_s$ (resp. $g_u$) the semisimple (resp. unipotent) part of $g$. If $G$ is
reductive we set $G_{\ad}=G/(\text{centre of }G^0)$. 

For a group $\G$ let $\fe(\G)$ be the set of conjugacy classes in $\G$. If $\G$ is a finite
group let $\Irr(\G)$ be the set of isomorphism classes of irreducible representations of 
$\G$ over $\bbq$.

For any finite set $I$ we denote by $\ZZ/2[I]$ the $\ZZ/2$-vector space with basis $I$. If
$I\ne\em$ we denote by $\ZZ/2[I]^\bul$ the codimension $1$ subspace of $\ZZ/2[I]$ 
consisting of vectors whose sum of coordinates is $0$; if $I=\em$ we set
$\ZZ/2[I]^\bul=\ZZ/2[I]=\{0\}$.

\head Contents\endhead
1. Preliminaries.

2. Semisimple classes in finite spin grpups.

3. Counting conjugacy classes in finite spin groups.

4. Counting data in the dual group.

5. Classification of irreducible representations.

\head 1. Preliminaries \endhead
\subhead 1.1\endsubhead
Let $\kk^\tri=\kk-\{0,1,-1\}$. Let $J=\{\l\in\kk;\l^2=-1\}\sub\kk^\tri$.
Define commuting permutations $\a,\b,\g$ of $\kk^\tri$ by $\a(\l)=\l\i$, $\b(\l)=-\l$,
$\g(\l)=\l^q$. We have $\a^2=\b^2=1$. For $e\in\ZZ/2$ we set $\g_e=\b^e\g$. 

Let $\O_\a$ (resp. $\O_{\g_e}$) be the set of orbits of $\a:\kk^\tri@>>>\kk^\tri$ (resp. 
$\g_e:\kk^\tri@>>>\kk^\tri$). Let $\O_{\a,\g_e}$ be the set of orbits of the group 
generated by $\a,\g_e$ acting on $\kk^\tri$. Let $\O_{\a,\b,\g}$ be the set of orbits of
the group generated by $\a,\b,\g$ acting on $\kk^\tri$. 

Let $\O'_{\a,\g_e}$ be the set of all $\co\in\O_{\a,\g_e}$ such that $\co$ is a single 
$\g_e$-orbit in $\kk^\tri$. Let $\O''_{\a,\g_e}=\O_{\a,\g_e}-\O'_{\a,\g_e}$; if 
$\co\in\O''_{\a,\g_e}$ then $\co$ is a union of two $\g_e$-orbits in $\kk^\tri$.

Let $X$ be an indeterminate. Let $\ca=\{\D\in\kk[X];\D\text{ monic },\D(0)\ne0\}$. We 
have $\ca=\sqc_{n\in\NN}\ca^n$ where $\ca^n=\{\D\in\ca;\deg(\D)=n\}$. Note that 
$\ca^0=\{1\}$. For $\D\in\ca$ we define a function $\kk^*@>>>\NN$, $\l\m n^\D_\l$ by
$$\D=\prod_{\l\in\kk^*}(X-\l)^{n^\D_\l}.$$ 

If $\D\in\ca$ and $\co\sub\kk^*$ is such that $\l\m n^\D_\l$ is constant on $\co$ then 
we set $n^\D_\co=n^\D_\l$ where $\l\in\co$.

For $\D\in\ca$ and $\x\in\{\a,\b,\g,\b\g\}$ we set 
${}^\x\D=\prod_{\l\in\kk^*}(X-\x(\l))^{n^\D_\l}$. 
$$\fA=\{\D\in\ca;{}^\a\D=\D,n^\D_1\in2\NN,n^\D_{-1}\in2\NN\},$$
$$\fA_0=\{\D\in\ca;{}^\a\D=\D,n^\D_1=0,n^\D_{-1}=0\}\sub\fA.$$

For $n\in\NN$ let $\fA^n=\fA\cap\ca^n$, $\fA^n_0=\fA_0\cap\ca^n$. For $n$ odd we have 
$\fA^n=\fA^n_0=\em$. Note that $\fA^n,\fA^n_0$ are stable under $\D\m{}^\b\D$, 
$\D\m{}^\g\D$. 

For $e\in\ZZ/2$ we set ${}^{\g_e}\fA^n=\{\D\in\fA^n;{}^{\g_e}\D=\D\}$,
${}^{\g_e}\fA^n_0=\{\D\in\fA^n_0;{}^{\g_e}\D=\D\}$.

\subhead 1.2\endsubhead
Let $e\in\ZZ/2$. For any $n\in\NN$ the condition that $\D\in\ca^n$ satisfies 
${}^{\g_e}\D=\D$ is that it is of the form $X^n+a_1X^{n-1}+\do+a_{n-1}X+a_n$ where 
$a_i\in\kk$ satisfy $(-1)^{ei}a_i^q=a_i$ (and $a_0\ne0$ if $n>0$). Hence the number of such
$\D$ is equal to $q^n-q^{n-1}$ (if $n>0$) and to $1$ if $n=0$. This number may also be 
computed in terms of the multiset of roots of $\D$ and then it appears as the number of 
solutions of the equation $\sum_{\co\in\O_{g_e}}|\co|n_\co+n_1+n_{-1}=n$ with 
$n_\co,n_1,n_{-1}$ in $\NN$ and, in the case where $e=1$, with $n_1=n_{-1}$. We see that
$$(1-qX)\i(1-X)=(1-qX)\i-X(1-qX)\i=1+\sum_{n\ge1}(q^n-q^{n-1})X^n$$
is equal to 
$$\prod_{\co\in\O_\g}(1-X^{|\co|})\i(1-X)^{-2}$$
if $e=0$ and to
$$\prod_{\co\in\O_{\g_1}}(1-X^{|\co|})\i(1-X^2)\i$$ if $e=1$.

\subhead 1.3\endsubhead
For any $n\in2\NN$, the condition that $\D\in\fA^n$ satisfies ${}^{\g_e}\D=\D$ is that it
is of the form $X^n+a_1X^{n-1}+\do+a_{n-1}X+a_n$ where $a_i\in\kk$ satisfy $a_i=a_{n-i}$   
for $i\in[1,n-1]$, $a_n=1$, $(-1)^{ei}a_i^q=a_i$. Hence the number of such $\D$ is equal to
$q^{n/2}$. This number may also be computed in terms of the multiset of roots of $\D$ and 
then it appears as the number of solutions of the equation
$$\sum_{\co\in\O_{\a,\g_e}}|\co|n_\co+n_1+n_{-1}=n$$
with $n_\co\in\NN$, $n_1,n_{-1}$ in $2\NN$ and (in the case where $e=1$) with $n_1=n_{-1}$.
We see that
$$(1-qX^2)\i=\sum_{n\in2\NN}q^{n/2}X^n$$
is equal to
$$\prod_{\co\in\O_{\a,\g}}(1-X^{|\co|})\i(1-X^2)^{-2}$$
if $e=0$ and to
$$\prod_{\co\in\O_{\a,\g_1}}(1-X^{|\co|})\i(1-X^4)\i$$
if $e=1$. 

\subhead 1.4\endsubhead
From 1.2 we have
$$\prod_{\co\in\O'_{\a,\g}}(1-X^{|\co|})\i\prod_{\co\in\O''_{\a,\g}}(1-X^{|\co|/2})^{-2}
=(1-qX)\i(1-X)^3,\tag a$$
$$\prod_{\co\in\O'_{\a,\g_1}}(1-X^{|\co|})\i\prod_{\co\in\O''_{\a,\g_1}}(1-X^{|\co|/2})^{-2}
=(1-qX)\i(1-X)(1-X^2).\tag b$$
From 1.3 we have
$$\prod_{\co\in\O'_{\a,\g}}(1-X^{|\co|})\i\prod_{\co\in\O''_{\a,\g}}(1-X^{|\co|})\i
=(1-qX^2)\i(1-X^2)^2,\tag c$$
$$\prod_{\co\in\O'_{\a,\g_1}}(1-X^{|\co|})\i\prod_{\co\in\O''_{\a,\g_1}}(1-X^{|\co|})\i
=(1-qX^2)\i(1-X^4).\tag d$$
In (a) we replace $X$ by $X^2$ and divide the resulting identity term by term with (c). We
obtain
$$\prod_{\co\in\O'_{\a,\g}}(1+X^{|\co|})\i\prod_{\co\in\O''_{\a,\g}}(1-X^{|\co|})\i=1-X^2.
\tag e$$
In (b) we replace $X$ by $X^2$ and divide the resulting identity term by term with (d). We
obtain
$$\prod_{\co\in\O'_{\a,\g_1}}(1+X^{|\co|})\i\prod_{\co\in\O''_{\a,\g_1}}(1-X^{|\co|})\i
=1-X^2.\tag f$$
Let $e\in\ZZ/2$.  For $\co\in\O_{\a,\g_e}$ we set ${}^ej(\co)=1$ if $\co\in\O'_{\a,\g_e}$
and ${}^ej(\co)=0$ if $\co\in\O''_{\a,\g_e}$. For $n\in2\NN,\D\in{}^{\g_e}\fA^n$ we set 
$${}^ej_\D=\sum_{\co\in\O_{\a,\g_e}}{}^ej(\co)n_\co^\D
=\sum_{\co\in\O'_{\a,\g_e}}n_\co^\D\in\ZZ/2.$$
When $e=0$ we shall also write $j(\co),j_\D$ instead of ${}^0j(\co),{}^0j_\D$.
We have ${}^ej_\D={}^ej_{{}^\b\D}$. From (e),(f) we see that for $n\in2\NN$:

(g) {\it the sum $\sum_{\D\in{}^{\g_e}\fA^n_0}(-1)^{{}^ej_\D}$ is $0$ if $n\ge4$, is $-1$
if $n=2$ and is $1$ if $n=0$.}
\nl
We show:

(h) {\it Let $e\in\ZZ/2$. Let $n\in2\NN,\D\in\fA^n$ be such that ${}^\b\D=\D={}^\g\D$.
Then ${}^ej_\D=n(q-(-1)^e)/4\in\ZZ/2$.}
\nl
We write $\O,\O'$ instead of $\O_{\a,\g_e}$, $\O'_{\a,\g_e}$. For $\l\in\kk^*$ we set
$n_\l=n^\D_\l$. For $\co\in\O$, $n_\co=n^\D_\co$ is defined. We must show that 
$$\sum_{\co\in\O'}n_\co=n(q-(-1)^e)/4\in\ZZ/2.$$
If $\co\in\O'$ then $\b(\co)\in\O'$ and $n_{\b(\co)}=n_{\co}$. Assume that $\co\in\O'$,
$\b(\co)=\co$. Let $\l\in\co$. Then $\b(\l)=\g_e^i(\l)$ for some $i\in[1,|\co|-1]$ and 
$\g_e^{2i}(\l)=\l$. Thus $|\co|$ divides $2i$. Hence $i=|\co|/2$ and $\b(\l)=\a(\l)$ that
is $-\l=\l\i$. We see that $\co\sub J$. (See 1.1.) Since $|\co|\ge2=|J|$ we see that 
$\co=J$. In particular, $J\in\O'$ so that $(-1)^e\l^q=\l\i$ for $\l\in J$ that is 
$q=-(-1)^e\mod4$. Thus the involution $\co\m\b(\co)$ of $\O'$ has $1$ fixed point $J$ if
$q=-(-1)^e\mod4$ and no fixed point if $q=(-1)^e\mod4$. We see that $\sum_{\co\in\O'}n_\co$
equals $n_J$ if $q=-(-1)^e\mod4$ and equals $0$ if $q=(-1)^e\mod4$. Hence (h) holds for 
$q=(-1)^e\mod4$. Now assume that $q=-(-1)^e\mod4$. We have $n=N'+N''+2n_1$ where 
$N'=\sum_{\l\in J}n_\l$, $N''=\sum_{\l\in\kk^\tri-J}n_\l$. Now $\kk^\tri-J$ is a disjoint 
union of four element sets $\{\l,\l\i,-\l,-\l\i\}$ on which $\l\m n_\l$ is constant. Hence
$N''\in4\NN$. Since $n_1\in2\NN$ we have $2n_1\in4\NN$. Hence $n=2n_J\mod4$ and 
$n_J=N/2\mod2$. This proves (h).

\subhead 1.5\endsubhead
Let $\co\in\O''_{\a,\g}$. If $\l\in\co$ then $\l^{q^{|\co|/2}}=l$. Hence 
$\l^{(q^{|\co|/2}-1)/2}\in\{1,-1\}$ and we can define $\e(\co)\in\ZZ/2$ by
$(-1)^{\e(\co)}=\l^{(q^{|\co|/2}-1)/2}$ (it is independent of the choice of $\l$ in $\co$).
Let $\co\in\O'_{\a,\g}$. If $\l\in\co$ then $\l^{q^{|\co|/2}}=\l\i$. Hence 
$\l^{(q^{|\co|/2}+1)/2}\in\{1,-1\}$ and we can define $\e(\co)\in\ZZ/2$ by
$(-1)^{\e(\co)}=\l^{(q^{|\co|/2}+1)/2}$ (it is independent of the choice of $\l$ in $\co$).
Thus $\e(\co)\in\ZZ/2$ is defined for any $\co\in\O_{\a,\g}$.

For any $n\in2\NN$ and $\D\in{}^\g\fA^n$ we set
$$\e_\D=(1/4)(q-1)n^\D_{-1}+\sum_{\co\in\O_{\a,\g}}\e(\co)n_\co^\D\in\ZZ/2.$$
We set

$\uu=(-1)^{(q-1)/2}$.
\nl
We now state the following result.

(a) {\it Let $x_{n;0}:=\sum_{\D\in{}^\g\fA^n_0;\e_D=0}(-1)^{j_\D}$. We have $x_{n;0}=0$ if 
$n\ge4$; $x_{n;0}=(1/2)(-1-\uu)$ if $n=2$; $x_{n;0}=1$ if $n=0$.}
\nl
The proof is given in 2.10. We now derive some consequences of (a). 

(b) {\it Let $x_{n;1}:=\sum_{\D\in{}^\g\fA^n_0;\e_D=1}(-1)^{j_\D}$. We have $x_{n;1}=0$ if 
$n\ge4$; $x_{n;1}=(1/2)(-1+\uu)$ if $n=2$; $x_{n;1}=0$ if $n=0$.}

(c) {\it Let $x_n:=\sum_{\D\in{}^\g\fA^N_0}(-1)^{j_\D+\e_\D}$. We have $x_n=0$ if $n\ge4$;
$x_n=-\uu$ if $n=2$; $x_n=1$ if $n=0$.}

(d) $\prod_{\co\in\O_{\a,\g}}(1-((-1)^{j(\co)+\e(\co)}X^{|\co|}))\i=1-\uu X^2$ .
\nl
Now (b) follows from (a) since, by 1.4(g), $x_{n;0}+x_{n;1}$ is equal to $0$ if $n\ge4$, to 
$-1$ if $n=2$ and to $1$ if $n=0$. We have $x_n=x_{n;0}-x_{n;1}$ hence (c) follows from
(a),(b). Clearly (d) follows from (c).

\subhead 1.6\endsubhead
For any $n\in2\NN$ the condition that $\D\in\fA^n$ satisfies $\D={}^\b\D={}^\g\D$ is that 
it is of the form $X^n+a_1X^{n-1}+\do+a_{n-1}X+a_n$ where $a_i\in\kk$ satisfy $a_i=0$ for 
odd $i$, $a_i=a_{n-i}$ for $i\in[1,n-1],a_n=1,a_i^q=a_i$. Hence the number of such $\D$ is
equal to $q^{n/4}$ if $n\in4\NN$ and to $q^{(n-2)/4}$ (if $n\in2+4\NN$). This number may 
also be computed in terms of the multiset of roots of $\D$ and then it appears as the 
number of solutions of $\sum_{\co\in\O_{\a,\b,\g}}|\co|n_\co+n_1+n_{-1}=n$ with 
$n_\co\in\NN$ and $n_1=n_{-1}\in2\NN$. We see that
$$\align&\prod_{\co\in\O_{\a,\b,\g}}(1-X^{|\co|})\i(1-X^4)\i=
\sum_{n\in4\NN}q^{n/4}X^n+\sum_{n\in2+4\NN}q^{(n-2)/4}X^n\\&=
(1-qX^4)\i+X^2(1-qX^4)\i=(1-qX^4)\i(1+X^2).\endalign$$
Hence
$$\prod_{\co\in\O_{\a,\b,\g};\co\ne J}(1-X^{|\co|})\i=(1-qX^4)\i(1-X^4)^2.\tag a$$

\subhead 1.7\endsubhead
Let $\cu$ be the set of all subsets $U$ of $\kk^\tri$ such that $\kk^\tri=U\sqc\a(U)$. 

Let $S$ be a set such that $|S|=2$. Let $n\in2\NN,\D\in\fA^n_0$. Let $S_\D$ be the set of 
all maps $\ps:\cu@>>>S$ such that

$\ps(U')=n/2+\sum_{\l\in U\cap U'}n^\D_\l+\ps(U)$ for any $U,U'$ in $\cu$.
\nl
We show:

(a) {\it for any $U\in\cu$, the map $S_\D@>>>S$, $\ps\m\ps(U)$ is a bijection.}
\nl
Let $n_\l=n^\D_\l$. It is enough to show that for $U,U',U''$ in $\cu$ we have
$$\sum_{\l\in U\cap U'}n_\l+\sum_{\l\in U'\cap U''}n_\l+\sum_{\l\in U''\cap U}n_\l
=\sum_{\l\in U}n_\l+\sum_{\l\in U'}n_\l+\sum_{\l\in U''}n_\l\mod2$$
or equivalently
$$\align&\sum_{\l\in(U\cap U')-U''}n_\l+\sum_{\l\in U''-(U\cup U')}n_\l+
\sum_{\l\in(U'\cap U'')-U}n_\l\\&+\sum_{\l\in U-(U'\cup U'')}n_\l+
\sum_{\l\in(U''\cap U)-U'}n_\l+\sum_{\l\in U'-(U''\cap U)}n_\l=0\mod2.\endalign$$
This follows from the fact that $\a$ defines bijections
$$\align&(U\cap U')-U''@>\si>>U''-(U\cup U'),(U'\cap U'')-U@>\si>>U-(U'\cup U''),\\&
(U''\cap U)-U'@>\si>>U'-(U''\cap U).\endalign$$
We see that 

$|S_\D|=|S|=2$. 

\subhead 1.8\endsubhead
Let $\cu^\bul=\{U\in\cu;U=\a\b(U)\}$. Note that $\cu^\bul\ne\em$.  

Let $n\in2\NN,\D\in\fA^n_0$. Note that $n^{{}^\b\D}_\l=n^\D_{\b(\l)}$ for any $\l$. To any 
$\ps\in S_\D$ and any $U\in\cu^\bul$ we associate an element ${}^\b\ps\in S_{{}^\b\D}$ by
the requirement that ${}^\b\ps(U)=n/2+\ps(U)$, see 1.7(a). We show that ${}^\b\ps$ is
independent of the choice of $U$. It is enough to show that if $U'\in\cu^\bul$ then 
$$\ps(U')=n/2+\sum_{\l\in U\cap U'}n^{{}^\b\D}_\l+\ps(U).$$
Since $\ps(U')=n/2+\sum_{\l\in U\cap U'}n^\D_\l+\ps(U)$, it is enough to show that 
$$\sum_{\l\in U\cap U'}n^{{}^\b\D}_\l=\sum_{\l\in U\cap U'}n^\D_\l.$$
The left hand side is 
$$\sum_{\l\in U\cap U'}n^\D_{\b(\l)}=\sum_{\l\in\b(U)\cap\b(U')}n^\D_\l=
\sum_{\l\in\a(U)\cap\a(U')}n^\D_\l=\sum_{\l\in U\cap U'}n^\D_{a(\l)}=
\sum_{\l\in U\cap U'}n^\D_\l,$$
as required.

The map $\ps\m{}^\b\ps$ is a bijection $S_\D@>>>S_{{}^\b\D}$.

\subhead 1.9\endsubhead
Let $e\in\ZZ/2$. Let $\cu^\bst_e$ be the set of all subsets $U$ of $\kk^\tri$ such that

(i) if $\co\in\O'_{\a,\g_e}$ then $U\cap\co=\{\g_e^i(\l);i=0,1,\do,|\co|/2-1\}$ for some
$\l\in\co$;

(ii) if $\co\in\O''_{\a,\g_e}$ then $\co\cap U$ is a single $\g_e$-orbit in $\kk^\tri$.
\nl
Note that in the setup of (i), $\l$ is uniquely determined by $U$; we set 
$\z_{\co,U}=\g_e\i(\l)$. We have $\em\ne\cu^\bst_e\sub\cu$. 

Let $n\in2\NN,\D\in\fA^n_0$. To any $\ps\in S_\D$ and any $U\in\cu^\bst_e$ we associate 
an element ${}^{\g_e}\ps\in S_{{}^{\g_e}\D}$ by the requirement that 
$${}^{\g_e}\ps(U)=\sum_{\co\in\O'_{\a,\g_e}}n^\D_{\z_{\co,U}}+\ps(U).\tag a$$
We show that ${}^{\g_e}\ps$ is independent of the choice of $U$. For $\l\in\kk^\tri$ we 
write $n_\l=n^\D_\l$, $m_\l=n^{{}^{\g_e}\D}_\l$. We have $m_{\g_e(\l)}=n_\l$. It is enough
to show that if $U'\in\cu^\bst_e$, then 
$$\sum_{\co\in\O'_{\a,\g_e}}n_{\z_{\co,U'}}+\ps(U')=
n/2+\sum_{\l\in U\cap U'}m_\l+\sum_{\co\in\O'_{\a,\g_e}}n_{\z_{\co,U}}+\ps(U).$$
Since $\ps(U')=n/2+\sum_{\l\in U\cap U'}n_\l+\ps(U)$, it is enough to show that 
$$\sum_{\co\in\O'_{\a,\g_e}}n_{\z_{\co,U'}}+\sum_{\l\in U\cap U'}n_\l=
\sum_{\l\in U\cap U'}m_\l+\sum_{\co\in\O'_{\a,\g_e}}n_{\z_{\co,U}}\mod2.$$
If $\co\in\O''_{\a,\g_e}$ then $U\cap U'\cap\co$ is $\g_e$-stable hence 
$$\sum_{\l\in U\cap U'\cap\co}n_\l=\sum_{\l\in U\cap U'\cap\co}m_\l.$$
Hence it is enough to show that for any $\co\in\O'_{\a,\g_e}$ we have
$$n_{\z_{\co,U'}}+\sum_{\l\in U\cap U'\cap\co}n_\l=
\sum_{\l\in\g_e\i(U\cap U'\cap\co)}n_\l+n_{\z_{\co,U}}\mod2.$$
We may identify $\co=\ZZ/(2r)$ so that the permutation $\g_e$ of $\co$ is $i\m i+1$. Then
$U\cap\co$ (resp. $U'\cap\co$) becomes the image of $\{a,a+1,\do,a+r-1\}$ (resp. 
$\{b,b+1,\do,b+r-1\}$) under $\ZZ@>>>\ZZ/(2r)$; $\z_{\co,U}$ (resp. $\z_{\co,U'}$) becomes
the image of $a-1$ (resp. $b-1$) under $\ZZ@>>>\ZZ/(2r)$ and the restriction of $\l\m n_\l$
to $\co$ becomes a function $n:\ZZ/(2r)@>>>\NN$ such that $n(i)=n(i+r)$ for all $i$. We can
assume that either $a\le b\le a+r-1\le b+r-1$ or $b\le a\le b+r-1\le a+r-1$. We can 
identify $U\cap U'$ with $\{b,b+1,\do,a+r-1\}$ (in the first case) and with 
$\{a,a+1,\do,b+r-1\}$ (in the second case). In the first case we must show that 
$$n_{b-1}+(n_b+n_{b+1}+\do+n_{a+r-1})=(n_{b-1}+n_b+\do+n_{a+r-2})+n_{a-1}\mod2$$
which is obvious since $n_{a-1}=n_{a+r-1}$. In the second case we must show that 
$$n_{b-1}+(n_a+n_{a+1}+\do+n_{b+r-1})=(n_{a-1}+n_a+\do+n_{b+r-2})+n_{a-1}\mod2$$
which is obvious since $n_{b-1}=n_{b+r-1}$. 

The map $\ps\m{}^{\g_e}\ps$ is a bijection $S_\D@>\si>>S_{{}^{\g_e}\D}$.

\head 2. Semisimple classes in finite spin grpups\endhead
\subhead 2.1\endsubhead
Let $V$ be a $\kk$-vector space of even dimension $N$ with a fixed nondegenerate 
quadratic form $Q:V@>>>\kk$. Let 
$$\align&O_Q=\{g\in GL(V);Q(g(x))=Q(x)\text{ for all }x\in V\},\\&
SO_Q=\{g\in O_Q;\det(g)=1\},\qua PSO_Q=SO_Q/\{1,-1\}.\endalign$$ 
For any subspace $V'$ of $V$ we set 

$V'{}^\pe=\{v\in V;Q(v+v')=Q(v)+Q(v')\text{ for all }v'\in V'\}$. 
\nl
Let $C(Q)$ be the Clifford algebra of $Q$. Recall that $C(Q)$ is an associative 
$\kk$-algebra with $1$ with a given $\kk$-linear imbedding $V@>>>C(Q)$ such that $v^2=Q(v)$
for any $v\in V$. Let $\tO_Q$ be the subgroup of the group of units of $C(Q)$ generated by
the elements $v\in V$ such that $Q(v)=1$. Note that $\tO_Q$ is a reductive algebraic group
(a closed subgroup of the group of units of $C(Q)$.)

If $v\in V$ and $g\in\tO_Q$ then $gvg\i$ (product in $C(Q)$) is in $V$ and 
$g\m[v\m gvg\i]$ is a homomorphism of algebraic groups 

$\k:\tO_Q@>>>O_Q$. 
\nl
For any 
$v\in V$ such that $Q(v)=1$ we have $\k(v)(v)=v$ and $\k(v)(v')=-v'$ if $v'\in(\kk v)^\pe$.
Let $\tSO_Q=\k\i(SO_Q)$ be the spin group of $Q$.
Now $\k$ restricts to a (surjective) homomorphism of algebraic groups 
$\tSO_Q@>>>SO_Q$ denoted again by $\k$; if $N\ge2$, its kernel has order $2$. For 
$N\ge2$ we define $\d\in\tSO_Q$ by $\k(\d)=1,\d\ne1$.

Let $\Xi=\{s\in SO_Q;s\text{ semisimple}\}$. For $s\in\Xi$ and $\l\in\kk^*$, let 

$V^s_\l=\{v\in V;sv=\l v\}$;
\nl
let 

$\D_s=\prod_{\l\in\kk^*}(X-\l)^{\dim V^s_\l}$.
\nl
We have $\D_s\in\fA^N$. We have $V=\op_{\l\in\kk^*}V^y_\l$ and $Q|_{V^s_1}$, 
$Q|_{V^s_{-1}}$ are nondegenerate. 

For $y\in SO_Q$ and $\l\in\{1,-1\}$ let
$$I^y_\l=\{[a];a\text{ odd, }y_u|_{V^{y_s}_\l}\text{ has some Jordan block of size }a\};$$
let 
$$I^y=\{[a];a\text{ odd, }y_u|_{V^{y_s}_1+V^{y_s}_{-1}}\text{ has some Jordan block of 
size }a\}.$$
Assuming that $N\ge2$
let $\tcs$ be the set of subspaces $E$ of $V$ such that $\dim E=N/2$ and $Q|_E=0$. Let
$\cs$ be the set of orbits of the $SO_Q$-action $y:E\m y(E)$ on $\tcs$. Let $E\m[[E]]$ be 
the obvious map $\tcs@>>>\cs$. Since $|\cs|=2$, the definition and results in 1.7 are
applicable to $S=\cs$. If $s\in\Xi$ satisfies $V^s_1=V^s_{-1}=0$ (so that $\D_s\in\fA^N_0$)
we define $\ps_s:\cu@>>>\cs$ by 

$\ps_s(U)=[[\op_{\l\in U}V^s_\l]]$.
\nl
We show that 
$$\ps_s\in\cs_{\D_s}.$$
(Notation of 1.7.) It is enough to show that for any $U,U'$ in $\cu$ we have
$$[[\op_{\l\in U'}V^s_\l]]=N/2+\sum_{\l\in U\cap U'}\dim V^s_\l+
[[\op_{\l\in U}V^s_\l]].$$
This follows from the following known result: 

(a) {\it if $E,E'\in\tcs$ then $[[E']]=\dim(E/(E\cap E'))+[[E]]$.}

\subhead 2.2\endsubhead
Assume that $N\ge2$.
For any $\D\in\fA^N$ we set $R_\D=\{s\in\Xi;\D_s=\D\}$. For any $\D\in\fA^N_0$ and any 
$\ps\in\cs_\D$ we set $R_\D^\ps=\{s\in R_\D;\ps_s=\ps\}$. Now (a),(b) below are easily
verified.

(a) {\it $R_\D$ ($\D\in\fA^N$) are exactly the $O_Q$-conjugacy classes in $\Xi$.}

(b) {\it $R_\D$ ($\D\in\fA^N-\fA^N_0$) and $R_\D^\ps$ ($\D\in\fA^N_0,\ps\in\cs_\D$) are 
exactly the $SO_Q$-conjugacy classes in $\Xi$.}

\subhead 2.3\endsubhead
Assume that $N\ge2$.
For any $g\in\tSO_Q$ (resp. $y\in SO_Q$) we denote by $Z(g)$ (resp. $Z(y)$) the 
centralizer of $g$ (resp. $y$) in $\tSO_Q$ (resp. in $SO_Q$); we also set 
$\tZ(g)=\{x\in\tSO_Q;xgx\i\in\{g,\d g\}\}=\k\i(Z(\k(g))$, a subgroup containing $Z(g)$.

We show:

(a) {\it Let $g\in\tSO_Q,y=\k(y)$. We have $xgx\i=\d g$ for some $x\in\tSO_Q$ if and only 
if $y$ satisfies conditions (i),(ii) below:

(i) $Z(y_s)$ is disconnected that is, $V^{y_s}_1\ne0$ and $V^{y_s}_{-1}\ne0$;

(ii) the $Z(y_s)^0$-conjugacy class of $y_u$ is also a $Z(y_s)$-conjugacy class that is, 
$I^y_1\ne\em$ and $I^y_{-1}\ne\em$.}
\nl
We have $\k\i(Z(y_s))=\tZ(g_s)$. Since $\k:\tSO_Q@>>>SO_Q$ is a double covering we have 
$\dim Z(g_s)=\dim\k(Z(g_s))$ and $\dim\tZ(g_s)=\dim Z(y_s)$ hence 
$\dim Z(y_s)^0=\dim\tZ(g_s)^0$. Since $|\tZ(g_s)/Z(g_s)|\le2$ we have 
$Z(g_s)^0=\tZ(g_s)^0$. Since $Z(g_s)=Z(g_s)^0$ we have $\k(Z(g_s))\sub Z(y_s)^0$ and 
$Z(g_s)=\tZ(g_s)^0$. Hence $\dim Z(g_s)=\dim Z(y_s)^0=\dim\k(Z(g_s))$. Thus 
$\k(Z(g_s))=Z(y_s)^0$. Since $\d\in Z(g_s)$ we have $Z(g_s)=\k\i(Z(y_s)^0)$. In particular,
$|Z(y_s)/Z(y_s)^0|=|\tZ(g_s)/Z(g_s)|\le2$.

Now assume that we can find $x\in\tSO_Q$ such that $xgx\i=\d g$. Then
$xg_sx\i=\d g_s,xg_ux\i=g_u$. Thus we have $x\n Z(g_s)$ so that $x\n\k\i(Z(y_s)^0)$ and
$\k(x)\n Z(y_s)^0$. Clearly, $\k(x)\in Z(y_s)$. Thus $Z(y_s)\ne Z(y_s)^0$ and $Z(y_s)$ is
disconnected. Also, if $z\in Z(y_s)$ we have either $z\in Z(y_s)^0$ or 
$z\in Z(y_s)^0\k(x)$ (note that $|Z(y_s)/Z(y_s)^0|\le2$); since $\k(x)y_u\k(x)\i=y_u$, in
both cases $zy_uz\i$ is in the $Z(y_s)^0$-conjugacy class of $y_u$. We see that $y$
satisfies (i),(ii).

Conversely, assume that $y$ satisfies (i),(ii). We have $|Z(y_s)/Z(y_s)^0|=2$ hence 
$|\tZ(g_s)/Z(g_s)|=2$. Then we can find $\x\in\tZ(g_s)$ such that $\x g_s\x\i=\d g_s$. We
have $\k(\x)\in Z(y_s)$ and since $y$ satisfies (ii) we have $\k(\x)y_u\k(\x)\i=zy_uz\i$ 
for some $z\in Z(y_s)^0$. We have $z=\k(\z)$ for some $\z\in Z(g_s)$ and 
$\k(\x)\k(g_u)\k(\x)\i=\k(\z)\k(g_u)\k(\z)\i$. Hence $\x g_u\x\i=\d^j\z g_u\z\i$ for some
$j\in\ZZ$. Since $\x g_u\x\i,\z g_u\z\i$ are unipotent and $\d^j$ is semisimple, central,
we must have $\d^j=1$ and $\x g_u\x\i=\z g_u\z\i$. Setting $\x'=\z\i\x$ we have 
$\x'g_u\x'{}\i=g_u$, $\x'g_s\x'{}\i=\d g_s$ hence $\x'g\x'{}\i=\d g$. This proves (a).

From (a) we deduce (b),(c),(d) below:

(b) {\it If $\D\in\fA^N$ and both $n^D_1,n^\D_{-1}$ are $\ne0$ then $\k\i(R_\D)$ is a 
single $\tSO_Q$-conjugacy class.}

(c) {\it If $\D\in\fA^N$ and exactly one of $n^D_1,n^\D_{-1}$ is $\ne0$ then $\k\i(R_\D)$
is a union of two $\tSO_Q$-conjugacy classes interchanged by multiplication by $\d$.}

(d) {\it If $\D\in\fA^N_0$ and $\ps\in\cs_\D$ then $\k\i(R_\D^\ps)$ is a union of two 
$\tSO_Q$-conjugacy classes interchanged by multiplication by $\d$.}
\nl
Indeed, the centralizer in $SO_Q$ of an element in $R_\D$ is disconnected in case (b) and 
connected in case (b); the centralizer in $SO_Q$ of an element in $R_\D^\ps$ in case (d) 
is connected.

\subhead 2.4\endsubhead
Let $\cf$ be the set of all group isomorphisms $F:V@>>>V$ such that $F(\l v)=\l^qF(v)$ for 
all $\l\in\kk,v\in V$ and $Q(F(v))=Q(v)^q$ for any $v\in V$. Let $F\in\cf$. Then $V^F$ is 
an $\FF_q$-vector space of dimension $N$ and $Q|_{V^F}$ is a nondegenerate quadratic form 
$V^F@>>>\FF_q$. We write $[V]_F$ for the image of $(V^F,Q|_{V^F})$ in $\Wi^0(\FF_q)=\ZZ/2$.
In the case where $N\ge2$ an equivalent definition of $[V]_F$ is that
$[[F(E)]]=[V]_F+[[E]]$ for any $E\in\tcs$. (In the case $N=0$ we have $[V]_F=0$.)

We have a partition $\cf=\cf_0\sqc\cf_1$ where for $e\in\ZZ/2$ we set 
$\cf_e=\{F'\in\cf;[V]_{F'}=e\}$.

Define $F:O_Q@>>>O_Q$ by $F(g)(F(v))=F(g(v))$ for $g\in O_Q,v\in V$; this restricts to a 
group isomorphism $SO_Q@>>>SO_Q$ denoted again by $F$ and this induces a group 
isomorphism $PSO_Q@>>>PSO_Q$ denoted again by $F$. There is a unique ring isomorphism 
$C(Q)@>>>C(Q)$ whose restriction to $V$ is $F:V@>>>V$ and whose restriction to 
$\kk\sub C(Q)$ is $\l\m\l^q$; we denote it again by $F$. This restricts to a group 
isomorphism $\tO_Q@>>>\tO_Q$ and to a group isomorphism $\tSO_Q@>>>\tSO_Q$ which are
denoted again by $F$. 

\subhead 2.5\endsubhead
Assume that $N\ge2$. Let $e\in\ZZ/2$. 
Let $F\in\cf$. Let $s\in\Xi$. Let $\D=\D_s$. Then $\D_{(-1)^eF(s)}={}^{\g_e}\D$. 
Now assume that $s$ satisfies $V^s_1=V^s_{-1}=0$.
Then $\ps_s\in\cs_\D$ and $\ps_{(-1)^eF(s)}\in{}^{\g_e}\D$,
${}^{\g_e}\ps_s\in\cs_{{}^{\g_e}\D}$ are defined. We show:
$$\ps_{(-1)^eF(s)}=[V]_F+{}^{\g_e}\ps_s.\tag a$$
It is enough to show that if $U\in\cu^\bst_e$ then (with notation of 1.9):
$$[[\op_{\l\in U}V^{(-1)^eF(s)}_\l]]
=[V]_F+\sum_{\co\in\O'_{\a,\g_e}}n^\D_{\z_{\co,U}}+[[\op_{\l\in U}V^s_\l]].\tag b$$
We have 
$$\op_{\l\in U}V^{(-1)^eF(s)}_\l=F(\op_{\l\in U}V^s_{\g_e\i(\l)})=
F(\op_{\l\in\g_e\i(U)}V^s_\l),$$
$$[[F(\op_{\l\in\g_e\i(U)}V^s_\l)]]=[V]_F+[[\op_{\l\in\g_e\i(U)}V^s_\l]],$$
$$[[\op_{\l\in\g_e\i(U)}V^s_\l]]=N/2+\dim(\op_{\l\in\g_e\i(U)\cap U}V^s_\l)
+[[\op_{\l\in U}V^s_\l]].$$
It is enough to show:
$$\sum_{\l\in\g_e\i(U)\cap U}n^\D_\l+\sum_{\co\in\O'_{\a,\g_e}}n^\D_{\z_{\co,U}}=N/2\mod2.
$$
From definitions we see that the left hand side is equal to $\sum_{\l\in U}n^\D_\l$ and 
this is clearly equal to $N/2$ as an integer, not only modulo $2$.

\subhead 2.6\endsubhead
Assume that $N\ge2$. Let $e\in\ZZ/2$. Let $F\in\cf$. Let $s\in\Xi$ be such that 
$(-1)^eF(s)=s$, $\D=\D_s$. Then $\D={}^{\g_e}\D$. Let $V'=V^s_1+V^s_{-1}$. Note that 
$Q|_{V'}$ is nondegenerate and $V'$ is $F$-stable. We show:
$${}^ej_\D=[V]_F-[V']_F\in\ZZ/2.\tag a$$
When $V=V'$ both sides are $0$. Hence we may assume that $V\ne V'$. Replacing if necessary
$V$ by $V'{}^\pe$, we may assume that $V'=0$ that is, $\D\in\fA^N_0$. Let $U\in\cu^\bst_e$.
In our case we may rewrite 2.5(b) as follows:
$$[[\op_{\l\in U}V^s_\l]]=[V]_F+{}^ej_\D+[[\op_{\l\in U}V^s_\l]].$$
We see that $[V]_F+{}^ej_\D=0$, as desired.

We show:

(b) {\it If in addition $\D={}^\b\D$ and $\D\in\fA^N_0$ then 
$[V]_F=N(q-(-1)^e)/4\in\ZZ/2$.}
\nl
This follows from (a) and 1.4(h).

\subhead 2.7\endsubhead
{\it In the remainder of this paper we assume that $N\ge2$.}

Let $F\in\cf$. For $\D\in\fA^N$ we have $F(R_\D)=R_{{}^\g\D}$. Hence $R_\D$ is $F$-stable
if and only if $\D\in{}^\g\fA^N$. 

If $\D\in{}^\g\fA^N$ and $n^D_1\ne0,n^\D_{-1}\ne0$ then $R_\D$ is a single ($F$-stable)
$SO_Q$-conjugacy class and for each $y$ in this class, $Z(y)$ has two connected components.
Hence in this case, $R_\D^F$ is a union of two $SO_Q^F$-conjugacy classes.

If $\D\in{}^\g\fA^N$ and exactly one of $n^D_1,n^\D_{-1}$ is nonzero then $R_\D$ is a 
single ($F$-stable) $SO_Q$-conjugacy class and for each $y$ in this class, $Z(y)$ is
connected. Hence in this case, $R_\D^F$ is a single $SO_Q^F$-conjugacy class.

If $\D\in{}^\g\fA^N_0$ then $R_\D$ is $F$-stable is a union of two $SO_Q$-conjugacy classes
$R_\D^\ps (\ps\in\cs_\D)$. We show:

(a) {\it Let $a=[V]_F+j_\D\in\ZZ/2$. For $\ps\in\cs_\D$ we have 
$F(R_\D^\ps)=R_\D^{a+\ps}$.}
\nl
Let $s\in R_\D^\ps$. We have $\D_s=\D,\ps_s=\ps$. We have $\D_{F(s)}={}^\g\D$ hence 
$\D_{F(s)}=\D$. By 2.5(a), 1.9(a) we have 

$\ps_{F(s)}(U)=[V]_F+{}^\g\ps_s(U)=[V]_F+j_\D+\ps_s(U)=a+\ps_s(U)$
\nl
for any $U\in\cu^\bst_0$. Thus $\ps_{F(s)}(U)=a+\ps_s(U)$ so that 
$\ps_{F(s)}=a+\ps_s=a+\ps$. We see that $F(s)\in R_\D^{a+\ps}$. This proves (a).

We see that if $a=0$ then $R_\D^\ps$ is $F$-stable. Since for any $y$ in this class, $Z(y)$
is connected we see that in this case $(R_\D^\ps)^F$ is a single $SO_Q^F$-conjugacy class.

\subhead 2.8\endsubhead
Let $F\in\cf$. Let $\fN_F:SO_Q^F@>>>\ZZ/2$ be the spinor norm. This is the group 
homomorphism defined by the requirement that $F(g)=\d^{\fN_F(y)}g$ for any $y\in SO_Q^F$,
$g\in\k\i(y)$. We have the following result.

(a) {\it Let $y\in SO_Q^F$. Setting $\D:=\D_{y_s}\in{}^\g\fA^N$, $e'=[V^{y_s}_{-1}]_F$, we 
have $(-1)^{\fN_F(y)}=(-1)^{e'+\e_\D}$.}
\nl
(Notation of 1.5.) Clearly any unipotent element of $SO_Q^F$ is in the kernel of $\fN_F$. 
Hence it is enough to prove (a) assuming in addition that $y=y_s$. If $V'$ is an 
$F$-stable, $y$-stable subspace of $V$ such that $Q|_{V'}$ is nondegenerate and such that 
$V'\ne0,V'\ne V$; from the definitions we see that, if (a) holds for $(V',y|_{V'})$ and for
$(V'{}^\pe,y|_{V'{}^\pe}$) then it also holds for $(V,y)$. Thus we may assume that there is
no subspace $V'$ as above. Then we are in one of the following four cases.

(i) $N=2$, $y=1$;

(ii) $N=2$, $y=-1$;

(iii) $V^y_1=V^y_{-1}=0$ and there is a unique $\co\in\O_{\a,\g}$ such that $n^\D_\co>0$;
moreover $n^\D_\co=1$ and $\co\in\O'_{\a,\g}$; 

(iv) $V^y_1=V^y_{-1}=0$ and there is a unique $\co\in\O_{\a,\g}$ such that $n^\D_\co>0$;
moreover $n^\D_\co=1$ and $\co\in\O''_{\a,\g}$.
\nl
In case (i) the result is trivial. 

Assume that we are in case (ii). Choose $\z\in\FF_q^*$ so that $\z^{(q-1)/2}=(-1)^{e'}$. We
can find a basis $u_1,u_2$ of $V^F$ such that $Q(au_1+a'u_2)=a^2-a'{}^2\z$ for 
$a,a'\in\kk^*$. Choose $\l\in\kk^*$ so that $\l^2=-\z$. We set $v=\l\i u_2$. We have 
$Q(v)=1$. Hence $g:=u_1v\in\k\i(-1)$ and $g=\l\i u_1u_2$, $F(g)=(\l\i)^qu_1u_2=bg$ where 

$b=(\l\i)^q\l=(\l^2)^{(1-q)/2}=(-\z)^{(1-q)/2}=(-1)^{(q-1)/2}(-1)^{e'}$. 
\nl
Hence (a) holds in this case.

Assume that we are in case (iii). Then $T:=Z(y)$ is an $F$-stable maximal torus of $SO_Q$
and $T^F$ is cyclic of order $q^{|\co|/2}+1$. It follows that $T':=\k\i(T)$ is an
$F$-stable maximal torus of $\tSO_Q$ and $T'{}^F$ is necessarily cyclic of order 
$q^{|\co|/2}+1$. Since $\d\in T'{}^F$ the map $T'{}^F@>>>T^F$ induced by $\k$ has kernel 
$\{1,\d\}$ hence is image is the unique subgroup of index $2$ of $T^F$. Hence $y\in T^F$ is
in the image of this map if and only if $y^{(q^{|\co|/2}+1)/2}=1$. Thus (a) holds in this 
case.

In case (iv) the proof is the same as in case (iii) provided that
we replace $q^{|\co|/2}+1$ by
$q^{|\co|/2}-1$. This completes the proof of (a).

\subhead 2.9\endsubhead
{\it In the remainder of this paper we assume that we have chosen $F_0\in\cf_0$ and
$F_1\in\cf_1$.} Let $e\in\ZZ/2$. For any $\D\in{}^\g\fA^N$ let $f_\D^e$ be the number
of semisimple conjugacy classes $\cc$ in $\tSO_Q^{F_e}$ such that $\k(\cc)\sub R_\D$; let 
$f_\D=f_D^0-f_\D^1$. 

We compute $f_\D$ in various cases. Set $n=n_1^\D$, $n'=n_{-1}^\D$.

(i) Assume that $n\ne0,n'\ne0$. In this case $\k\i(R_\D)$ is a single ($F_e$-stable)
semisimple conjugacy class in $\tSO_Q$ hence $(\k\i(R_\D))^{F_e}$ is a single conjugacy 
class in $\tSO_Q^{F_e}$ (we use the fact that the centralizer of a semisimple element in 
$\tSO_Q$ is connected). Hence we have $f_\D^e=1$. We see that $f_\D=0$.

(ii) Assume that $n\ne0,n'=0$. In this case $R_\D^{F_e}$ is a single conjugacy class $\boc$
in $SO_Q^{F_e}$. By 2.8(a), the spinor norm is equal to $(-1)^{\e_\D}$ on $\boc$.
Hence if $\e_\D=1$ then $f_\D^0=f_\D^1=0$ and $f_\D=0$, while if $\e_\D=0$ we have 
$\k\i(\boc)^{F_e}\ne\em$. Since by 2.3(c), for any $g\in\k\i(\boc)$, $\d g$ is not 
conjugate to $g$ under $\tSO_Q$ (hence under $\tSO_Q^{F_e}$),  we see that when $\e_\D=0$,
$\k\i(\boc)^{F_e}$ is a union of two conjugacy classes in $\tSO_Q^{F_e}$. Thus if 
$\e_\D=0$, we have $f_\D^0=f_\D^1=2$ and $f_\D=0$.

(iii) Assume that $n=0,n'\ne0$. In this case $R_\D^{F_e}$ is a single conjugacy class 
$\boc$ in $SO_Q^{F_e}$. By 2.8(a), the spinor norm is equal to $(-1)^{e'+\e_\D}$ on $\boc$
where $e'=[V^s_{-1}]_F$ with $s\in\boc$. By 2.6(a) we have $e'=e+j_\D$. Thus the spinor
norm is equal to $(-1)^{e+j_\D+\e_\D}$ on $\boc$. As in case (ii) we see that $f_\D^e$
equals $0$ if $e+j_\D+\e_\D=1$ and $f_\D^e$ equals $2$ if $e+j_\D+\e_\D=0$. Thus 
$f_\D=2(-1)^{j_\D+\e_\D}$.

(iv) Assume that $n=n'=0$. If $e\ne j_\D$ the $f_\D^e=0$, see 2.7(a). Now assume that
$e=j_\D$. By 2.7(a), for any $\ps\in\cs_\D$, $(R_\D^\ps)^{F_e}$ is a single conjugacy class
$\boc$ in $SO_Q^{F_e}$. By 2.8(a), the spinor norm is equal to $(-1)^{\e_\D}$ on $\boc$. 
As in case (ii) we see that $f_\D^e$ equals $0$ if $\e_\D=1$ and $f_\D^e$ equals $4$ if 
$\e_\D=0$. (We have $4$ instead of $2$ since there are two $\ps$'s.) We see that 
$f_\D=4(-1)^{j_\D}$ if $\e_\D=0$ and $f_\D=0$ if $\e_\D=1$.

\subhead 2.10. Proof of 1.5(a)\endsubhead
Recall from 1.5 the notation $x_{n;i}=\sum_{\D\in{}^\g\fA^N_0;\e_\D=i}(-1)^{j_\D}$ for
$i=0,1$. Let $f=\sum_{\D\in{}^\g\fA^N}f_\D$. We compute $f_\D$ in two different ways. On 
the one hand we have $f=f_0-f_1$ where $f_e$ is the number of semisimple conjugacy classes
in $\tSO_Q^{F_e}$. Since $\tSO_Q$ is semisimple simply-connected of rank
$N/2$ for $N\ge4$, we have $f_e=q^{N/2}$ hence $f=0$; for $N=2$ we have $f_e=q-(-1)^e$,
hence $f=-2$. 

On the other hand from 2.9 we see that $f_\D=0$ if $n^\D_1\ne0$. If $n^\D_1=0$ we can write
uniquely $\D=(X+1)^{n'}\D'$ where $n'\in2\NN$ and $\D'\in{}^\g\fA^{N-n'}_0$; note that 
$j_\D=j_{\D'},\e_\D=\e_{\D'}+n'(q-1)/4$. Using this and 2.9 we obtain   
$$f=\sum_{\D\in{}^\g\fA^N_0;\e_\D=0}4(-1)^{j_\D}+\sum_{n'\in2\NN;0<n'\le N}
\sum_{\D'\in{}^\g\fA^{N-n'}_0}2(-1)^{j_{\D'}+\e_{\D'}}\uu^{n'/2}$$
that is
$$0=4x_{N;0}+2\sum_{N'\in2\NN;N'<N}(x_{N';0}-x_{N';1})\uu^{(N-N')/2}$$
for $N\ge4$ and $-2=4x_{2;0}+2(x_{0;0}-x_{0,1})\uu$. Recall that by 1.4(g), 
$x_{n;0}+x_{n;1}$ is equal to $0$ if $n\ge4$ and to $-1$ if $n=2$. Also from the 
definitions we have $x_{0;0}=1$, $x_{0;1}=0$. It follows that
$$0=4x_{N;0}+4\sum_{N'\in2\NN;2<N'<N}x_{N';0}\uu^{(N-N')/2}+
2(2x_{2;0}+1)\uu^{(N-2)/2}+2\uu^{N/2}$$
for $N\ge4$ and $x_{2,0}=(-2-2\uu)/4$. Hence for $N\ge4$ we have

$2(2x_{2;0}+1)\uu^{(N-2)/2}+2\uu^{N/2}=0$
\nl
and 
$x_{N;0}+\sum_{N'\in2\NN;2<N'<N}x_{N';0}\uu^{(N-N')/2}=0$.
\nl
From this we see by induction on $N$ that $x_{N;0}=0$ if $N\ge4$. This proves 1.5(a).

\head 3. Counting conjugacy classes in finite spin groups\endhead
\subhead 3.1\endsubhead
{\it Here and in 3.2, 3.3 we fix $g\in\tSO_Q$ and we set $y=\k(g)$.}

We assume that there exist $V',V'',f$ where $V',V''$ are $y$-stable subspaces of $V$ both 
of dimension $a$ (odd) such that $V''\sub V'{}^\pe$, $Q|_{V'},Q|_{V''}$ are nondegenerate 
and $f:V'@>>>V''$ is an isometry such that $fy(v)=yf(v)$ for any $v\in V'$. We show:

(a) $\d\in Z(g)^0$.
\nl
Let $v'_1,\do,v'_a$ be a basis of $V'$ such that $Q(\sum_ia_iv'_i)=\sum_ia_i^2$ for any
$a_i\in\kk$. Let $v''_1,\do,v''_a$ be the basis of $V''$ given by $v''_i=f(v'_i)$. For any 
$\l,\mu$ in $\kk$ such that $\l^2+\mu^2=1$, the vectors $w_i=\l v'_i+\mu v''_i (i\in[1,a])$
form a basis for the subspace $V_{\l,\mu}$ they generate such that
$Q(\sum_ia_iw_i)=\sum_ia_i^2$ for any $a_i\in\kk$; moreover, $Q|_{V_{\l,\mu}}$ is 
nondegenerate and $V_{\l,\mu}$ is $y$-stable. Let $g_{\l,\mu}=(u_1u_2)^2\in\tSO_Q$ where 
$u_1=v'_1v'_2\do v'_a\in\tO_Q$, $u_2=w_1w_2\do w_a\in\tO_Q$. Now $\k(u_1)$ is $1$ on $V'$
and is $-1$ on $(V')^\pe$; hence it commutes with $y$. Similarly, $\k(u_2)$ is $1$ on 
$V_{\l,\mu}$ and is $-1$ on $(V_{\l,\mu})^\pe$; hence it commutes with $y$. Thus,
$\k(g_{\l,\mu})\in Z(y)$ and $g_{\l,\mu}\in\k\i(Z(y))=\tZ(g)$. The map 
$\{(\l,\mu)\in\kk\T\kk;\l^2+\mu^2=1\}@>>>\tZ(g)$, $(\l,\mu)\m g_{\l,\mu}$, has as image an
irreducible subset of $\tZ(g)$. Since $g_{1,0}=1$, this image must be contained in 
$\tZ(g)^0=Z(g)^0$. Since $a$ is odd we have $g_{0,1}=\d$. Hence (a) holds.

\subhead 3.2\endsubhead
We assume that there exist two $y$-stable subspaces $V',V''$ of $V$ both of dimension $a$ 
(odd) such that $V'\cap V''=0$, $Q|_{V'}=0$, $Q|_{V''}=0$, $Q|_{V'+V''}$ is nondegenerate.
We show:

(a) $\d\in Z(g)^0$.
\nl
Let $v_1,\do,v_a$ be a basis of $V'$. Let $v'_1,\do,v'_a$ be the basis of $V''$ such that
$Q(v_i+v'_j)=0$ for $i\ne j$ and $Q(v_i+v'_i)=1$ for all $i$. For any $c\in\kk^*$ we have 
$Q(cv_i+c\i v'_i)=1$. For $c\in\kk^*$ and $i=1,\do,a$ we set
$$g_{i,c}=(c\i-c)v_iv'_i+c=(c-c\i)v'_iv_i+c\i=(v_i+v'_i)(cv_i+c\i v'_i)\in\tSO_Q.$$
Let $y_{i,c}=\k(g_{i,c})$. We have $y_{i,c}(v_i)=c^{-2}v_i$, $y_{i,c}(v'_i)=c^2v'_i$ and 
$y_{i,c}(v)=v$ for any $v\in(\kk v_i+\kk v'_i)^\pe$. Moreover $g_{1,c},g_{2,c},\do,g_{a,c}$
commute with each other in $\tSO_Q$. Their product in $\tSO_Q$ is denoted by by $g_c$. Let
$y_c=\k(g_c)$. We have $y_c(v_i)=c^{-2}v_i$, $y_c(v'_i)=c^2v'_i$ for all $i$ and $y_c(v)=v$
for any $v\in(V'+V'')^\pe$. Thus $y_c\in Z(y)$ and $g_c\in\tZ(g)$. The map 
$\kk^*@>>>\tZ(g)$, $c\m g_c$ has as image an irreducible subset of $\tZ(g)$. Since $g_1=1$,
this image must be contained in $\tZ(g)^0=Z(g)^0$. Since $g_{-1}=(-1)^a=-1$ we see that 
(a) holds.

\subhead 3.3\endsubhead
We show:

(a) {\it If for some odd $a$, the number of Jordan blocks of size $a$ of $y_u:V@>>>V$ is 
$>1$ then $\d\in Z(g)^0$.}

(b) {\it If for any odd $a$, the number of Jordan blocks of size $a$ of $y_u:V@>>>V$ is
$\le1$ then $\d\n Z(g)^0$.}
\nl
For $a$ odd and $t\in\kk^*$, let $f_{a,t}$ be the number of Jordan blocks of size $a$ of
$y_u|_{V^{y_s}_t}$. Assume that the hypothesis of (a) holds. Then for some odd $a$ we have
$\sum_tf_{a,t}>1$. If $f_{a,t}>0$ for some $t\in\kk^\tri$ then we also have $f_{a,t\i}>0$.
We see that 3.2(a) is applicable, with $V'\sub V^{y_s}_t$ and $V''\sub V^{y_s}_{t\i}$, so 
that $\d\in Z(g)^0$. Next we assume that $f_{a,t}=0$ for any $t\in\kk^\tri$. Then 
$f_{a,1}+f_{a,-1}>1$. We see that 3.1(a) is applicable, with $V',V''$ subspaces of 
$V^{y_s}_1+V^{y_s}_{-1}$, so that $\d\in Z(g)^0$. This proves (a).

Assume now that the hypothesis of (b) holds. Then by an argument in \cite{\ICC, 14.3} we
have $\d\n Z(g_u)^0$. Since $Z(g)\sub Z(g_u)$ we have $Z(g)^0\sub Z(g_u)^0$ and we see that
$\d\n Z(g)^0$. This proves (b).

\subhead 3.4\endsubhead
Let $F\in\cf$. Let $C$ be a conjugacy class in $\tSO_Q$ such that $F(C)=C$. Then 
$C^F\ne\em$. We set $M_{C,F}=\{\cc\in\fe(\tSO_Q^F);\cc\sub C\}$. As it is known, one can 
describe $M_{C,F}$ in terms of the action of $F$ on $\ov{Z(g)}$ (where $g\in C^F$) as 
follows. Consider the equivalence relation $\si$ on $\ov{Z(g)}$ given by $h\si h'$ if 
$h'=h_1\i hF(h_1)$ for some $h_1\in\ov{Z(g)}$. We define a map $Z(g)@>>>M_{C,F}$ by 

$z\m(\tSO_Q^F-\text{conjugacy class of }\x g\x\i)$
\nl
where $\x\in\tSO_Q$ is such that $z=\x\i F(\x)$. This map is well defined and (by standard
arguments using Lang's theorem) it induces a bijection 

(a) $\ov{Z(g)}/\si@>>>M_{C,F}$.

\subhead 3.5\endsubhead
In the setup of 3.4 assume in addition that $\d C=C$. Then $\cc\m\d\cc$ defines an 
involution $\d_1$ of $M_{C,F}$. Now $h\m h_1\i hF(h_1)$ (with $h\in Z(g)/Z(g)^0$, 
$h_1\in\tZ(g)/Z(g)^0-Z(g)/Z(g)^0$) induces an involution $\d_0$ of $\ov{Z(g)}/\si$. We
have a commutative diagram
$$\CD
\ov{Z(g)}/\si@>>>M_{C,F}\\
@V\d_0VV                 @V\d_1VV\\
\ov{Z(g)}/\si@>>>M_{C,F}\endCD$$
where the horizontal maps are as in 3.4(a). Indeed, let $z\in Z(g)$ and let $\x\in\tSO_Q$ 
be such that $z=\x\i F(\x)$. Let $z'\in\tZ(g)-Z(g)$. Then $z'{}\i zF(z')\in Z(g)$ and 
$z'{}\i zF(z')=z'{}\i\x\i F(\x)F(z')$. It is enough to show that  
$(\x z')g(\x z')\i=\d\x g\x\i$. This follows from $z'gz'{}\i=\d g$. 

\subhead 3.6\endsubhead
Let $g\in\tSO_Q$ and let $y=\k(g)\in SO_Q$. We show:

(a) {\it the image of the homomorphism $Z(g)@>>>Z(y)$ induced by $\k$ is 
$Z(y_s)^0\cap Z(y)$; moreover the index of $Z(y_s)^0\cap Z(y)$ in $Z(y)$ is equal to the 
index of $Z(g)$ in $\tZ(g)$.}
\nl
Note that $\k(\tZ(g))=Z(y)$. Hence if $\tZ(g)=Z(g)$, (a) is clear (in this case
$Z(y_s)^0=Z(y_s)$ so that $Z(y_s)^0\cap Z(y)=Z(y)$). Now assume that $\tZ(g)\ne Z(g)$. We 
have $Z(g)\sub Z(g_s)=Z(g_s)^0$ hence 
$\k(Z(g))\sub\k(Z(g_s)^0)\sub Z(y_s)^0$. Also, $\k(Z(g))\sub Z(g)$. Hence 
$\k(Z(g))\sub Z(y_s)^0\cap Z(y)$. Since $|\tZ(g)/Z(g)|=2$ and $\ker\k\sub Z(g)$ we have 
$|\k(\tZ(g))/\k(Z(g))|=2$ that is, $|Z(y)/\k(Z(g))|=2$. It is enough to show that 
$Z(y)\not\sub Z(y_s)^0$. We can find $x\in\tZ(g)-Z(g)$. Then $xgx\i=\d g$ and 
$xg_sx\i=\d g_s$. Since $x\n Z(g_s)=\k\i(Z(y_s)^0)$, we have $\k(x)\n Z(y_s)^0$. On the 
other hand, $\k(x)\in Z(y)$. Thus $Z(y)\not\sub Z(y_s)^0$ and (a) is proved.

From (a) we see that we have an exact sequence

$1@>>>\{1,\d\}@>>>Z(g)@>\k>>Z(y_s)^0\cap Z(y)@>>>1$.
\nl
We see that we must be in one of the following two cases:

(i) $\d\in Z(g)^0$; then $\ov{Z(g)}=\ov{Z(y_s)^0\cap Z(y)}$;

(ii) $\d\n Z(g)^0$; then we have an exact sequence

$1@>>>\{1,\d\}@>>>\ov{Z(g)}@>>>\ov{Z(y_s)^0\cap Z(y)}@>>>1$.

\subhead 3.7\endsubhead
Let $F\in\cf$. Let $g\in\tSO_Q^F$ be such that $\d\n Z(g)^0$ and let $y=\k(g)\in SO_Q^F$. 
Using 3.3 we see that $I^y=I^y_1\cup I^y_{-1}$ is a disjoint union and that can find 
$y$-stable, $F$-stable subspaces $V_a ([a]\in I^y)$, $V'$ of $V$ which form a direct sum 
decomposition of $V$ such that 

any summand is contained in the ${}^\pe$ of any other summand;

$V_a\sub V^{y_s}_1$ if $[a]\in I^y_1$, $V_a\sub V^{y_s}_{-1}$ if $[a]\in I^y_{-1}$;

$\dim V_a=a$ and $y_u|_{V_a}$ has a single Jordan block if $[a]\in I^y$;

$y_u|_{V'}$ has no Jordan block of odd size.
\nl
Since $\dim V^{y_s}_1,\dim V^{y_s}_{-1}$ are even, we see that $|I^y_1|,|I^y_{-1}|$ are 
even.

For any $[a]\in I^y$ we define $r_a\in O_Q$ by $r_a(v)=v$ for $v\in V_a$, $r_a(v)=-v$ for
$v\in V_a^\pe$; let $\tir_a\in\tO_Q$ be one of the two elements in $\k\i(r_a)$. We have
$r_ay=yr_a$, $F(r_a)=r_a$. Using $(-1)^{aa'}=-1$ we see that 

(a) {\it if $[a],[a']$ in $I^y$ are distinct then $\tir_a\tir_{a'}=\d\tir_{a'}\tir a$.}
\nl
We identify $\ZZ/2[I^y]$ with the subgroup of $O_Q$ generated by $\{r_a;[a]\in I^y\}$ in
such a way that $[a]\in I^y$ corresponds to $r_a$. Then 
$E:=\ZZ/2[I^y_1]^\bul\op\ZZ/2[I^y_{-1}]^\bul$ is identified with a subgroup of 
$Z(y_s)^0\cap Z(y)$. Note that $E$ contains exactly one element in each connected component
of $Z(y_s)^0\cap Z(y)$. From 3.6 we see that $\tE:=\k\i(E)$ (a subgroup of $Z(g)$) contains
exactly one element in each connected component of $Z(g)$.

For $f,f'$ in $\ZZ/2[I^y]$ we define $(f,f')\in\ZZ/2$ by $\tf\tf'=\d^{(f,f')}\tf'\tf$ where
$\tf\in\k\i(f)$, $\tf'\in\k\i(f')$ are in $\tO_Q$; note that $(f,f')$ is well defined since
$\ZZ/2[I^y]$ is abelian (it is clearly independent of the choices). Now 
$(,):\ZZ/2[I^y]\T\ZZ/2[I^y]@>>>\ZZ/2$ is a symplectic bilinear form and by (a) we have 
$([a],[a'])=1$ for $[a]\ne[a']$. Let $\sum_{[a]\in I^y}x_a[a]$ be in the radical of $(,)$.
(Here $x_a\in\ZZ/2$.) Then $\sum_{[a]\in I^y;[a]\ne[a']}x_a=0$ for any $[a']\in I^y$. Hence
$x_{[a']}=\sum_{[a]}x_a$ is independent of $[a']$; we denote it by $x$. We have 
$x=|I^y|x=0$ since $|I^y|$ is even. Hence $x_{[a]}=0$ for all $[a]$. We see that $(,)$ is 
nondegenerate on $\ZZ/2[I^y]$.

Let $E_0$ be the subspace of $\ZZ/2[I^y]$ spanned by the elements 
$f_1=\sum_{[a]\in I^y_1}[a]$, $f_{-1}=\sum_{[a]\in I^y_{-1}}[a]$. We have 
$(f_1,f_{-1})=|I^y_1||I^y_{-1}|=0$ in $\ZZ/2$ hence $(,)|_{E_0}=0$. We have 
$E=\{f\in\ZZ/2[I^y];(f,f_1)=(f,f_{-1})=0\}$. Hence $E_0$ is exactly the radical of 
$(,)|_E$. 

Let $E'=\ZZ/2[I^y]^\bul$, $\tE'=\k\i(\tE)$. 
Note that $E\sub E'\sub Z(y)$, $E'$ contains exactly one element 
in each connected component of $Z(y)$ and $\tE'\sub\tZ(g)$ contains exactly one element in
each connected component of $\tZ(g)$.

Since $F(r_a)=r_a$ for every $[a]\in I^y$, the subgroup $\ZZ/2[I^y]$ of $O_Q$ is contained
in $O_Q^F$. Hence $\k\i(\ZZ/2[I^y])$ is an $F$-stable subgroup of $\tO_Q$ and for any 
$f\in\ZZ/2[I^y]$ we have $F(\tf)=\d^{c(f)}\tf$ (where $\tf\in\k\i(f)$) and $c(f)\in\ZZ/2$.
(Note that $c(f)$ is independent of the choice of $\tf$.) Clearly, $f\m c(f)$ is a 
homomorphism $\ZZ/2[I^y]@>>>\ZZ/2$.

Now $h':h\m h'{}\i hF(h')$ defines a $\tE'$-action on $\tE$ and, by restriction, a
$\tE$-action on $\tE$. Let $\cx$ be the set of orbits of the $\tE$-action. We show:

(b) {\it $|\cx|$ equals $|E|$ if $c(f')\ne0$ for some $f'\in E_0$ and equals $|E|+|E_0|$ if
$c|_{E_0}=0$.}
\nl
For $f,f'$ in $\ZZ/2[I^y]$ we have $\tf'{}\i\tf F(\tf')=\d^{(f,f')+c(f')}\tf$ where 
$\tf\in\k\i(f)$, $\tf'\in\k\i(f')$. Hence $\tE@>>>E$ induces a surjective map $\cx@>>>E$
whose fibre at $f\in E$ has size $2$ if $(f,f')+c(f')=0$ for all $f'\in E$ and has size $1$
otherwise. We see that $\cx=|E|+|A|$ where

$A=\{f\in E;(f,f')+c(f')=0\text{ for any }f'\in E\}$. 
\nl
Assume first that $c(f_1)\ne0$. Then for any $f\in E$ we have $(f,f_1)+c(f_1)\ne0$; thus 
$A=\em$. Similarly, if $c(f_{-1})\ne0$ we have $A=\em$. Next assume that 
$c(f_1)=c(f_{-1})=0$. Then $c|_{E_0}=0$. Hence we can find $\uf\in E$ such that 
$c(f')=(\uf,f')$ for any $f'\in E$. We have 
$A=\{f\in E;(f+\uf,f')=0\text{ for any }f'\in E\}=\uf+E_0$. Thus, $|A|=|E_0|$. This proves
(b).

In the case where $I^y_1\ne\em,I^y_{-1}\ne\em,$ we have a partition $\cx=\cx'\sqc\cx''$ 
where $\cx'$ (resp. $\cx''$) is the set of $\tE$-orbits on $\tE$ which are (resp. are not)
$\tE'$-orbits. In this case we have the following refinement of (b):

(c) {\it We have $|\cx'|=|E|$. Moreover $|\cx''|$ equals $0$ if $c(f')\ne0$ for some 
$f'\in E_0$ and equals $|E_0|=4$ if $c|_{E_0}=0$.}
\nl
From the first sentence in the proof of (b) we see that $|\cx'|$ is the number of 
$\tE$-orbits of size $2$ on $\tE$ plus the number of fixed points of the $\tE'$-action on 
$\tE$. Hence
$$\align&|\cx'|=|\{f\in E;(f,f')+c(f')\ne0\text{ for some }f'\in E\}|\\&
+2|\{f\in E;(f,f')+c(f')=0\text{ for any }f'\in E'\}|.\endalign$$
To prove that $|\cx'|$ is as in (c) it is enough to show that $|A|=2|\tA|$ where $A$ is as 
in the proof of (b) and $\tA=\{f\in E;(f,f')+c(f')=0\text{ for any }f'\in E'\}$.

Assume first that $c(f')\ne0$ for some $f'\in E_0$. Then for any $f\in E$ we have 
$(f,f')+c(f')\ne0$; thus $A=\tA=\em$ and (c) holds. Now assume that $c|_{E_0}=0$. Then we 
can find $\uf\in E$ such that $c(f')=(\uf,f')$ for any $f'\in E$. As in the proof of (b) we
have $|A|=|E_0|=4$. Since $c|_L=0$ there exists $\hf\in E'$ such that $c(f')=(\hf,f')$ for
any $f'\in E'$. For $f'\in E$ we have $(\hf,f')=(\uf,f')$ hence $\uf\in\hf+E_0$. Thus we 
may assume that $\uf=\hf$. Note that $\tA\sub A$. The condition that $f\in\hf+E_0$ is in 
$\tA$ is that for any $f'\in E'$ we have $(f,f')+c(f')=0$ that is $(f+\hf,f')=0$ or that 
$f+\hf\in L$. We see that $\tA=\hf+L$ so that $|\tA|=2$ and $|A|-2|\tA|=4-2\T2=0$. This 
proves (c).

Returning to the general case we note that for $\l\in\{1,-1\}$ we have

$c(f_\l)=[W_\l]_F+\dim(W_\l)(q-1)/4$
\nl
where $W_\l=\sum_{[a]\in I^y_\l}V_a$. (We apply 2.8(a) to $f_\l\in SO_Q$ instead of $y$.)
We set $W'_\l=V^{y_s}_\l\cap V'$. Then $V^{y_s}_\l=W_\l\op W'_\l$ and both summands are
$F$-stable and $y$-stable; moreover $y_u:W'_\l@>>>W'_\l$ is an orthogonal transformation
with no Jordan blocks of odd size. It follows that $\dim W'_\l\in4\NN$ and $[W'_\l]_F=0$ so
that $[W_\l]_F=[V^{y_s}_\l]_F$ and $\dim W_\l=\dim V^{y_s}_\l\mod 4$. We see that

$c(f_\l)=[V^{y_s}_\l]_F+\dim(V^{y_s}_\l)(q-1)/4$.
\nl
Since $y_u$ has no odd Jordan blocks on $V^{y_s}_{\l'}$ with $\l'\in\kk^\tri$ we see that
$n^\D_{\l'}\in2\NN$ for any $\l'\in\kk^\tri$ (where $\D=\D_{y_s}$). Hence 
$\e_\D=\dim(V^{y_s}_{-1})(q-1)/4$.
Since $y=\k(g)$ we have $\fN_F(y)=0$. Using 2.8(a) we deduce that
$[V^{y_s}_{-1}]_F+\e_\D=0$ hence $[V^{y_s}_{-1}]_F+\dim(V^{y_s}_{-1})(q-1)/4=0$.
Thus $c(f_{-1})=0$. It follows that 

$c(f_1)=c(f_1)+c(f_{-1})=
[V^{y_s}_1+V^{y_s}_{-1}]_F+\dim(V^{y_s}_1+\dim V^{y_s}_{-1})(q-1)/4$.
\nl
By 2.6(a) we have $[V^{y_s}_1+V^{y_s}_{-1}]_F=[V]_F+j_\D$. Using $n^\D_{\l'}\in2\NN$ for 
$\l'\in\kk^\tri$ we see that $j_\D=0$ and $\dim(V^{y_s}_1+\dim V^{y_s}_{-1})=\dim(V)\mod4$.
It follows that

$c(f_1)=[V]_F+N(q-1)/4$.
\nl
Now let $C$ be the $\tSO_Q$-conjugacy class of $g$.
We have canonically $M_{C,F}\lra\ov{Z(g)}/\si\lra\cx$ (see 3.4(a)); moreover, if
$\d C=C$ then the fixed point set of $\d_1:M_{C,F}@>>>M_{C,F}$ corresponds to the fixed 
point set of $\d_0:\ov{Z(g)}/\si@>>>\ov{Z(g)}/\si$ (see 3.5(a)) and this is in bijection 
with the subset $\cx'$ of $\cx$. We can summarize the results above as follows.

(d) {\it Assume that $I^y_1=I^y_{-1}=\em$. Then $|M_{C,F}|=2$.}

(e) {\it Assume that $I^y_1\ne\em, I^y_{-1}=\em$. Then}

     $|M_{C,F}|=2^{|I^y_1|-1}+2$ if $[V]_F+N(q-1)/4=0\in\ZZ/2$,

     $|M_{C,F}|=2^{|I^y_1|-1}$ if $[V]_F+N(q-1)/4=1\in\ZZ/2$.

(f) {\it Assume that $I^y_1=\em, I^y_{-1}\ne\em$. Then} $|M_{C,F}|=2^{|I^y_1|-1}+2$.

(g)  {\it Assume that $I^y_1\ne\em, I^y_{-1}\ne\em$. Then 

     $|M_{C,F}|=2^{|I^y_1|+|I^y_{-1}|-2}+4$ if $[V]_F+N(q-1)/4=0\in\ZZ/2$,

     $|M_{C,F}|=2^{|I^y_1|+|I^y_{-1}|-2}$ if $[V]_F+N(q-1)/4=1\in\ZZ/2$.
\nl
Moreover, the number of $\cc\in M_{C,F}$ such that $\d\cc=\cc$ is equal to
$2^{|I^y_1|+|I^y_{-1}|-2}$.}

\subhead 3.8\endsubhead
Let $F\in\cf$. Let $g\in\tSO_Q^F$ be such that $\d\in Z(g)^0$ and let $y=\k(g)\in SO_Q^F$. 
From 3.6(i) we see that $\ov{Z(g)}=\ov{Z(y_s)^0\cap Z(y)}$ and this may be identified with
$\ZZ/2[I^y_1]^\bul\op\ZZ/2[I^y_{-1}]^\bul$ in such a way that the action of $F$ on 
$\ov{Z(g)}$ becomes the trivial action on $\ZZ/2[I^y_1]^\bul\op\ZZ/2[I^y_{-1}]^\bul$. Hence
if $C$ is the $\tSO_Q$-conjugacy class of $g$ we have canonically 
$M_{C,F}\lra\ov{Z(g)}/\si\lra\ZZ/2[I^y_1]^\bul\op\ZZ/2[I^y_{-1}]^\bul$ (see 3.4(a)). If in
addition we have $\d C=C$ then, by 3.5, the involution $\d_1$ of $M(C,F)$ becomes the
involution of $\ov{Z(y_s)^0\cap Z(y)}$ given by conjugation by an element of $\ov{Z(y)}$
which is not in $\ov{Z(y_s)^0\cap Z(y)}$ and this last involution is trivial since
$\ov{Z(y)}$ is commutative.

We can summarize the results above as follows.

(a) {\it Assume that $I^y_1=I^y_{-1}=\em$. Then $|M_{C,F}|=1$.}

(b) {\it Assume that exactly one of $I^y_1,I^y_{-1}$ is $\ne\em$. Then}
     $|M_{C,F}|=2^{|I^y_1|+|I^y_{-1}|-1}$.

(c)  {\it Assume that $I^y_1\ne\em, I^y_{-1}\ne\em$. Then 
     $|M_{C,F}|=2^{|I^y_1|+|I^y_{-1}|-2}$.
\nl
Moreover, the number of $\cc\in M_{C,F}$ such that $\d\cc=\cc$ is equal to
$2^{|I^y_1|+|I^y_{-1}|-2}$.}

\subhead 3.9\endsubhead
For any $n\in\NN$ let $P_n$ be the set of partitions of $n$ and let $\p(n)=|P_n|$. Thus we
have $\p(0)=1$. We set $\p(n)=0$ when $n\in\QQ-\NN$. For $n\in\NN$ let $\tP_n$ be the set 
of partition of $n$ into even parts. We have $|\tP_n|=\p(n/2)$.
Let $T_N$ be the set of all $SO_Q$-conjugacy classes of unipotent elements $y\in SO_Q$ such
that $I^y\ne\em$. Let $T_N^0$ be the set of all $SO_Q$-conjugacy classes of unipotent 
elements $y\in SO_Q$ such that $I^y=\em$. Let $\tT_N$ be the set of all $SO_Q$-conjugacy 
classes of unipotent elements $y\in SO_Q$ such that $I^y\ne\em$ and such that for any odd
$a$, $y$ has at most one Jordan block of size $a$. Now $T^0_N$ has a natural fixed point 
free involution in which an orbit and its image are contained in the same $O_Q$-conjugacy
class; let $\uT^0_N$ be the set of orbits of this involution. 
For $n,n'$ in $2\NN+2$ we note that $\ZZ/2$ acts on $T^0_n\T T^0_{n'}$ (by the involution
above acting simultaneously on both factors); let $\un{T^0_n\T T^0_{n'}}$ be the set of 
orbits for this $\ZZ/2$-action. We define $\uT^0_0$ to be a set with one element. 

We have $|\uT^0_n|=\p(n/4)$ for all $n\in2\NN$. For $n\in2\NN+2$ we set $\tit_n=|\tT_n|$
and $\t_n=\sum_y2^{|I^y|-1}$ where $y$ runs over a set of representatives for the conjugacy
classes in $T_n$. We also set $\tit_0=0,\t_0=0$.

\subhead 3.10\endsubhead
Let $e\in\ZZ/2$. We fix $\D\in{}^\g\fA^N$. We set 
$n_1=n_1^\D,n_{-1}=n_{-1}^\D$. For $\co\in\O_{\a,\g}$ we set $n_\co=n^\D_\co$. 
Let $A_e^\D=\{g\in\tSO_Q^{F_e};\k(g)_s\in R_\D\}$.
Let $\uA_e^\D$ be the set of $\tSO_Q^{F_e}$-conjugacy classes contained in $A_e^\D$.
Let ${}'\uA_e^\D=\{\cc'\in\uA_e^\D;\d\cc'=\cc'\}$,
${}''\uA_e^\D=\{\cc'\in\uA_e^\D;\d\cc'\ne\cc'\}$. 

\subhead 3.11\endsubhead
In the setup of 3.10, assume that $n_1\ne0$, $n_{-1}\ne0$. 
From the results of 2.9, 3.7, 3.8 we see that

$|{}''\uA_e^\D|=h_1+h_2+h_3+h_4+h_5+h_6+h_7$ where 

$h_1=2\prod_{\co\in\O_{\a,\g}}|P_{n_\co}||\uT^0_{n_1}|\t_{n_{-1}}$
if $e=j_\D+e_\D$ and $h_1=0$ otherwise;

$h_2=4\prod_{\co\in\O_{\a,\g}}|\tP_{n_\co}||\uT^0_{n_1}|\tit_{n_{-1}}$
if $e=(n_1+n_{-1})(q-1)/4$ and $h_2=0$ otherwise;

$h_3=2\prod_{\co\in\O_{\a,\g}}|P_{n_\co}|\t_{n_1}|\uT^0_{n_{-1}}|$
if $\e_\D=0$ and $h_3=0$ otherwise;

$h_4=4\prod_{\co\in\O_{\a,\g}}|\tP_{n_\co}|\tit_{n_1}|\uT^0_{n_{-1}}|$
if $e=(n_1+n_{-1})(q-1)/4$ and $h_4=0$ otherwise;

$h_5=2\prod_{\co\in\O_{\a,\g}}|P_{n_\co}||\un{T^0_{n_1}\T T^0_{n_{-1}}}|$
if $\e_D=0$, $e=j_\D$ and $h_5=0$ otherwise;

$h_6=2\prod_{\co\in\O_{\a,\g}}|\tP_{n_\co}||\un{T^0_{n_1}\T T^0_{n_{-1}}}|$
if $e=0$ and $h_6=0$ otherwise;

$h_7=4\prod_{\co\in\O_{\a,\g}}|\tP_{n_\co}|\tit_{n_1}\tit_{n_{-1}}$
if $e=(n_1+n_{-1})(q-1)/4$ and $h_7=0$, otherwise.

It follows that

$$\align&|{}''\uA_0^\D|-|{}''\uA_1^\D|=
2\prod_{\co\in\O_{\a,\g}}\p(n_\co)\p(n_1/4)\t_{n_{-1}}(-1)^{j_\D+e_\D}    \\&
+4\d_{\e_\D=0}\prod_{\co\in\O_{\a,\g}}\p(n_\co)\p(n_1/4)\p(n_{-1}/4)(-1)^{j_\D+e_\D}  \\&
+4\prod_{\co\in\O_{\a,\g}}\p(n_\co/2)\p(n_1/4)\tit_{n_{-1}}(-1)^{(n_1+n_{-1})(q-1)/4} \\&
+4\prod_{\co\in\O_{\a,\g}}\p(n_\co/2)\tit_{n_1}\p(n_{-1}/4)(-1)^{(n_1+n_{-1})(q-1)/4} \\&
+4\prod_{\co\in\O_{\a,\g}}\p(n_\co/2)\p(n_1/4)\p(n_{-1}/4)(-1)^{(n_1+n_{-1})(q-1)/4}  \\&
+4\prod_{\co\in\O_{\a,\g}}\p(n_\co/2)\tit_{n_1}\tit_{n_{-1}}(-1)^{(n_1+n_{-1})(q-1)/4}.
\endalign$$

\subhead 3.12\endsubhead
In the setup of 3.10, assume that $n_1\ne0$, $n_{-1}=0$. 
From the results of 2.9, 3.7, 3.8 we see that 

$|{}''\uA_e^\D|=h_1+h_2+h_3+h_4$ where 

$h_1=2\prod_{\co\in\O_{\a,\g}}|P_{n_\co}||T^0_{n_1}|$
if $\e_\D=0,e=j_\D+e_\D$ and $h_1=0$ otherwise;

$h_2=2\prod_{\co\in\O_{\a,\g}}|\tP_{n_\co}||T^0_{n_1}|$
if $\e_\D=0,e=(n_1+n_{-1})(q-1)/4$ and $h_2=0$ otherwise;

$h_3=2\prod_{\co\in\O_{\a,\g}}|P_{n_\co}|\t_{n_1}$
if $\e_\D=0$ and $h_2=0$ otherwise;

$h_4=4\prod_{\co\in\O_{\a,\g}}|P_{n_\co}|\tit_{n_1}$
if $\e_\D=0,e=(n_1+n_{-1})(q-1)/4$ and $h_4=0$ otherwise.

It follows that

$$\align&|{}''\uA_0^\D|-|{}''\uA_1^\D|=
4\d_{\e_\D=0}\prod_{\co\in\O_{\a,\g}}\p(n_\co)\p(n_1/4)(-1)^{j_\D+e_\D}    \\&
+4\d_{\e_\D=0}
\prod_{\co\in\O_{\a,\g}}\p(n_\co/2)(\p(n_1/4)+\tit_{n_1})(-1)^{(n_1+n_{-1})(q-1)/4}.
\endalign$$

\subhead 3.13\endsubhead
In the setup of 3.10, assume that $n_1=0$, $n_{-1}\ne0$. 
From the results of 2.9, 3.7, 3.8 we see that 

$|{}''\uA_e^\D|=h_1+h_2+h_3+h_4$ where 

$h_1=2\prod_{\co\in\O_{\a,\g}}|P_{n_\co}||T^0_{n_{-1}}|$
if $\e_\D=0,e=j_\D+e_\D$ and $h_1=0$ otherwise;

$h_2=2\prod_{\co\in\O_{\a,\g}}|\tP_{n_\co}||T^0_{n_{-1}}|$
if $e=(n_1+n_{-1})(q-1)/4$ and $h_2=0$ otherwise;

$h_3=2\prod_{\co\in\O_{\a,\g}}|P_{n_\co}|\t_{n_{-1}}$
if $e=j_\D+e_\D$ and $h_3=0$ otherwise;

$h_4=4\prod_{\co\in\O_{\a,\g}}|\tP_{n_\co}|\tit_{n_{-1}}$
if $e=(n_1+n_{-1})(q-1)/4$ and $h_4=0$ otherwise;

It follows that

$$\align&|{}''\uA_0^\D|-|{}''\uA_1^\D|=
2\prod_{\co\in\O_{\a,\g}}\p(n_\co)\t_{n_{-1}}(-1)^{j_\D+e_\D}    \\&
+4\d_{\e_\D=0}\prod_{\co\in\O_{\a,\g}}\p(n_\co)\p(n_{-1}/4)(-1)^{j_\D+e_\D}    \\&
+4\prod_{\co\in\O_{\a,\g}}\p(n_\co/2)(\p(n_{-1}/4)+\tit_{n_{-1}})
(-1)^{(n_1+n_{-1})(q-1)/4}.\endalign$$

\subhead 3.14\endsubhead
In the setup of 3.10, assume that $n_1=n_{-1}=0$. From the results of 2.9, 3.7, 3.8 we see
that 

$|{}''\uA_e^\D|=h_1+h_2$ where 

$h_1=4\prod_{\co\in\O_{\a,\g}}|P_{n_\co}|$
if $\e_\D=0,e=j_\D+e_\D$ and $h_1=0$ otherwise;

$h_2=4\prod_{\co\in\O_{\a,\g}}|\tP_{n_\co}|$
if $\e_\D=0,e=(n_1+n_{-1})(q-1)/4$ and $h_2=0$ otherwise.

It follows that

$$\align&|{}''\uA_0^\D|-|{}''\uA_1^\D|=
4\d_{\e_\D=0}\prod_{\co\in\O_{\a,\g}}\p(n_\co)(-1)^{j_\D+e_\D}    \\&
+4\d_{\e_\D=0}\prod_{\co\in\O_{\a,\g}}\p(n_\co/2)(-1)^{(n_1+n_{-1})(q-1)/4}.\endalign$$

\subhead 3.15\endsubhead
We show:

(a) {\it If $C\in\fe(\tSO_Q)$ and $\d C=C$ then for any $\x\in\tO_Q$ we have $\x C\x\i=C$.}
\nl
Let $C'=\k(C)$. We have $C=\k\i(C')$. Let $\x\in\tO_Q$ and let $x=\k(\x)\in O_Q$. We must 
show that $\x C\x\i=C$. It is enough to show that $xC'x\i=C'$. Let $y\in C'$. By 2.3(a) 
$y$ satisfies (i),(ii) in
2.3(a). Hence we can find an orthogonal reflection $r:V@>>>V$ such that $r=1$ on 
$(V^{y_s}_1)^\pe$. It is enough to show that $ryr\i$ is $SO_Q$-conjugate to $y$. Now 
$(ryr\i)|_{V^{y_s}_1}$ and $y|_{V^{y_s}_1}$ are unipotent orthogonal transformations which 
are conjugate under $O_{Q|_{V^{y_s}_1}}$; moreover they have some Jordan block of odd size
hence we can find $u_1\in SO_{Q|_{V^{y_s}_1}}$ such that 
$(ryr\i)|_{V^{y_s}_1}=u_1(y|_{V^{y_s}_1})u_1\i$. Define $u\in SO_Q$ by $u(v)=u_1(v)$ for 
$v\in V^{y_s}_1$ and $u(v)=v$ for $v\in(V^{y_s}_1)^\pe$. Then $ryr\i=uyu\i$. This 
proves (a).

We show:

(b) {\it If $C$ is as in (a) then we have $F_0(C)=C$ if and only if $F_1(C)=C$.}
\nl
Let $y:=F_0\i F_1:V@>>>V$. Now $y$ is a vector space isomorphism preserving $Q$ hence 
$y\in O_Q$. We have $F_1(v)=F_0(y(v))$ for all $v\in\VV$. It follows that
$F_1(g)=F_0(\ty g\ty\i)$ for any $g\in\tO_Q$ where $\ty\in\tO_Q$ is such that $\k(\ty)=y$.
Hence the condition that $F_1(C)=C$ is equivalent to the condition that $F_0(\ty C\ty\i)=C$
and this is equivalent to the condition that $F(C)=C$, since $\ty C\ty\i=C$, by (a). This 
proves (b).

Let $C$ be as in (a) and such that $F_0(C)=C$ (or equivalently $F_1(C)=C$). By the results 
in 3.7, 3.8 we have

(c) $|M(C,F_0)|=|M(C,F_1)|=|\ov{\k(Z(g))}|$
\nl
for any $g\in C$. (Note that $|\ov{\k(Z(g))}|$ does not depend on the Frobenius map.)
Conbining this with (b) we see that in the setup of 3.10 we have
$$|{}'\uA_0^\D|=|{}'\uA_1^\D|.\tag d$$

\subhead 3.16\endsubhead
Let
$$\Ps(X)=\sum_{n\in\NN}\p(n)X^n=\prod_{k\ge1}(1-X^k)\i,$$
$$\La(X)=\sum_{n\in2\NN}\t_nX^n,\qua\ti\La(X)=\sum_{n\in2\NN}(\tit_n+\p(n/4))X^n.$$
In the identity
$$1+\sum_{n\in\NN+1}(2x'_n+(1/2)x''_n)X^n=\Ps(X^2)^2\sum_{m\in\ZZ}X^{m^2}\tag a$$
in \cite{\ICC, p.249} we have $x'_n=\t_n$ if $n\in2\NN$ and $x''_n=2\p(n/4)$ for all $n$.
Adding the identity (a) with the one obtained from it by the substitution $X\m-X$ we 
therefore obtain
$$\sum_{n\in2\NN}(4\t_n+2\p(n/4))X^n=2\Ps(X^2)^2\sum_{m\in2\ZZ}X^{m^2}$$
that is,
$$\sum_{n\in2\NN}(2\t_n+\p(n/4))X^n=\Ps(X^2)^2\sum_{m'\in\ZZ}X^{4m'{}^2}.$$
Using this and Jacobi's identity
$$\sum_{j\in\ZZ}X^{j^2}=\Ps(-X)^{-2}\Ps(X^2)\tag b$$
with $X$ replaced by $X^4$ we obtain
$$2\La(X)+\Ps(X^4)=\Ps(X^2)^2\Ps(-X^4)^{-2}\Ps(X^8)\tag c$$
In the identity
$$\sum_{n\in\NN}\sum_{t\in\ZZ}|\cp_{n,t}|X^n=\Ps(X^4)\Ps(X^2)\Ps(-X)\i\tag d$$
in \cite{\ICC, (14.5.2)} (with $Z_1=X$, $Z_2=1$) we have
$\sum_{t\in\ZZ}|\cp_{n,t}|=\tit_n+\p(n/4)$ if $n$ is even. Adding the identity (d) with the
one obtained from it by the substitution $X\m-X$ we therefore obtain
$$2\sum_{n\in2\NN}(\tit_n+\p(n/4))X^n=\Ps(X^4)\Ps(X^2)\Ps(-X)\i+\Ps(X^4)\Ps(X^2)\Ps(X)\i.$$
Thus
$$2\ti\La(X)=\Ps(X^4)\Ps(X^2)(\Ps(-X)\i+\Ps(X)\i).\tag e.$$
From the computation
$$\align&\prod_{k\ge1}(1-(-X)^k)\i
=\prod_{k\ge1;\text{ odd}}(1+X^k)\i\prod_{k\ge1;\text{ even}}(1-X^k)\i\\&
=\prod_{k\ge1}(1+X^k)\i\prod_{k\ge1;\text{ even}}(1+X^k)\prod_{k\ge1;\text{ even}}(1-X^k)\i
\\&=\prod_{k\ge1}(1-X^{2k})\i\prod_{k\ge1}(1-X^k)\prod_{k\ge1;\text{ even}}(1-X^{2k})\\&
\T\prod_{k\ge1;\text{ even}}(1-X^k)\i\prod_{k\ge1;\text{ even}}(1-X^k)\i\endalign$$
we see that
$$\Ps(-X)\Ps(X)=\Ps(X^2)^3\Ps(X^4)\i.\tag f$$

\subhead 3.17\endsubhead
We set 
$$\align&\a_N=|\Irr(\tSO_Q^{F_0})|-|\Irr(\tSO_Q^{F_1})|=
|\fe(\tSO_Q^{F_0})|-|\fe(\tSO_Q^{F_1})|\\&=\sum_{\D\in{}^\g\fA^N}(|\uA_0^\D|-|\uA_1^\D|)=
\sum_{\D\in{}^\g\fA^N}(|{}''\uA_0^\D|-|{}''\uA_1^\D|).\tag a\endalign$$
(The last equality follows from 3.15(d).) We set $\a_0=8$. From 3.11-3.14 we see that
$$\sum_{n\in2\NN}\a_nX^n=\fa+\fb+\fc$$
where
$$\align\fa&=2\sum_{\D\in{}^\g\fA}\prod_{\co\in\O_{\a,\g}}\p(n^\D_\co)
(-1)^{n^\D_\co j(\co)+n^\D_\co\e(\co)}X^{|\co|n^\D_\co}\\&
\T\p(n^\D_1/4)X^{n^\D_1}\t_{n^\D_{-1}}(-1)^{n^\D_{-1}(q-1)/4}X^{n^\D_{-1}},\endalign$$
$$\fb=4\sum_{\D\in{}^\g\fA,\e_D=0}\prod_{\co\in\O_{\a,\g}}\p(n^\D_\co)
(-1)^{n^\D_\co j(\co)}X^{|\co|n^\D_\co}
\T\p(n^\D_1/4)X^{n^\D_1}\p(n^\D_{-1}/4)X^{n^\D_{-1}},$$
$$\align&\fc=4\sum_{\D\in{}^\g\fA}\prod_{\co\in\O_{\a,\g}}\p(n^\D_\co/2)X^{|\co|n^\D_\co}\\&
\T(\tit_{n^\D_1}+\p(n^\D_1/4))X^{n^\D_1}(\tit_{n^\D_{-1}}+\p(n^\D_{-1}))
X^{n^\D_{-1}}(-1)^{(n^\D_1+n^\D_{-1})(q-1)/4}.\endalign$$
We set
$$\Ps_q(X)=\prod_{k\ge1}(1-qX^k)\i.\tag b$$
Recall that $\uu=(-1)^{(q-1)/2}$. Let $w$ be such that $w^2=\uu$. We have
$$\align&\fa=2(\sum_{n\in2\NN}\p(n/4)X^n)(\sum_{n\in2\NN}\t_n(wX)^n)
\prod_{\co\in\O_{\a,\g}}\{\sum_{n\in\NN}\p(n)((-1)^{j(\co)+\e(\co)}X^{|\co|})^n\}\\&
=2\Ps(X^4)\La(wX)
\prod_{\co\in\O_{\a,\g}}\prod_{k\ge1}(1-((-1)^{j(\co)+\e(\co)}X^{|\co|})^k)\i.\endalign$$
Thus
$$\align&2\i\Ps(X^4)\i\La(wX)\i\fa=
\prod_{k\ge1}\prod_{\co\in\O_{\a,\g}}(1-((-1)^{j(\co)+\e(\co)}X^{|\co|})^k)\i\\&
=\prod_{k\ge1;k\text{ odd}}\prod_{\co\in\O_{\a,\g}}(1-(-1)^{j(\co)+\e(\co)}(X^k)^{|\co|})\i
\prod_{k\ge2;k\text{ even}}\prod_{\co\in\O_{\a,\g}}(1-(X^k)^{|\co|})\i.\endalign$$
The products over $\co$ can be evaluated using 1.3 and 1.5(d). We obtain
$$\align&2\i\Ps(X^4)\i\La(wX)\i\fa=\prod_{k\ge1;k\text{ odd}}(1-\uu X^{2k})
\prod_{k'\ge1}(1-qX^{4k'})\i\prod_{k'\ge1}(1-X^{4k'})^2\\&
=\prod_{k\ge1,k\text{ odd}}(1-(wX)^{2k})\Ps_q(X^4)\Ps(X^4)^{-2}\\&
=\prod_{k\ge1}(1-(wX)^{2k})\prod_{k\ge2;k\text{ even}}(1-(wX)^{2k})\i
\Ps_q(X^4)\Ps(X^4)^{-2}\\&
=\Ps((wX)^2)\i\Ps(X^4)\Ps_q(X^4)\Ps(X^4)^{-2}=\Ps_q(X^4)\Ps(X^4)\i\Ps((wX)^2)\i.\endalign$$
Using this and 3.16(c) with $X$ replaced by $wX$ we obtain
$$\align&\fa=2\Ps_q(X^4)\Ps((wX)^2)\i\La(wX)\\&
=\Ps_q(X^4)\Ps((wX)^2)\Ps(-X^4)^{-2}\Ps(X^8)-\Ps_q(X^4)\Ps((wX)^2)\i\Ps(X^4).\endalign$$
We have $\fb=\fb'+\fb''$ where
$$\fb'=2\sum_{\D\in{}^\g\fA,\e_D=0}\prod_{\co\in\O_{\a,\g}}\p(n^\D_\co)
(-1)^{n^\D_\co j(\co)}X^{|\co|n^\D_\co}\p(n^\D_1/4)X^{n^\D_1}\p(n^\D_{-1}/4)X^{n^\D_{-1}},$$
$$\align\fb''&=2\sum_{\D\in{}^\g\fA,\e_D=0}\prod_{\co\in\O_{\a,\g}}\p(n^\D_\co)
(-1)^{n^\D_\co j(\co)}X^{|\co|n^\D_\co}\\&\p(n^\D_1/4)X^{n^\D_1}\p(n^\D_{-1}/4)X^{n^\D_{-1}}
(-1)^{\e_\D}.\endalign$$
We have
$$\align&\fb'=2(\sum_{n\in2\NN}\p(n/4)X^n)^2
\prod_{\co\in\O_{\a,\g}}\sum_{n\in\NN}\p(n)((-1)^{j(\co)}X^{|\co|})^n\\&
=2\Ps(X^4)^2\prod_{\co\in\O_{\a,\g}}\prod_{k\ge1}(1-((-1)^{j(\co)}X^{|\co|})^k)\i\\&
=2\Ps(X^4)^2\prod\Sb k\ge1;\\k\text{ odd}\endSb
\prod_{\co\in\O_{\a,\g}}(1-(-1)^{j(\co)}(X^k)^{|\co|})\i
\prod\Sb k\ge2;\\k\text{ even}\endSb\prod_{\co\in\O_{\a,\g}}(1-(X^k)^{|\co|})\i.\endalign$$
The products over $\co$ can be evaluated using 1.3 and 1.4(e). We obtain
$$\align&\fb'=2\Ps(X^4)^2\prod_{k\ge1;k\text{ odd}}(1-X^{2k})
\prod_{k\ge2;k\text{ even}}\{(1-qX^{2k})\i(1-X^{2k})^2\}\\&
=2\Ps(X^4)^2\prod_{k\ge1}(1-X^{2k})\prod_{k\ge2;k\text{ even}}(1-X^{2k})\i
\Ps_q(X^4)\Ps(X^4)^{-2}\\&=2\Ps(X^4)^2\Ps(X^2)\i\Ps(X^4)\Ps_q(X^4)\Ps(X^4)^{-2}
=2\Ps_q(X^4)\Ps(X^4)\Ps(X^2)\i.\endalign$$
We have
$$\align&\fb''=2\sum_{\D\in{}^\g\fA,\e_D=0}\prod_{\co\in\O_{\a,\g}}\p(n^\D_\co)
(-1)^{n^\D_\co(j(\co)+\e(\co))}\\&
\T X^{|\co|n^\D_\co}\p(n^\D_1/4)X^{n^\D_1}(-1)^{n^D_{-1}(q-1)/4}\p(n^\D_{-1}/4)X^{n^\D_{-1}}
\\&=2(\sum_{n\in2\NN}\p(n/4)X^n)^2
\prod_{\co\in\O_{\a,\g}}\sum_{n\in\NN}\p(n)((-1)^{j(\co)+\e(\co)}X^{|\co|})^n\\&
=2\Ps(X^4)^2\prod_{\co\in\O_{\a,\g}}\prod_{k\ge1}(1-((-1)^{j(\co)+\e_\co}X^{|\co|})^k)\i.
\endalign$$
As in the computation of $\fa$ the product over $\co$ is equal to
$$\Ps_q(X^4)\Ps(X^4)\i\Ps((wX)^2)\i.$$ Hence
$$\fb''=2\Ps_q(X^4)\Ps(X^4)\Ps((wX)^2)\i.$$
We have
$$\align&\fc=4\prod_{\co\in\O_{\a,\g}}\sum_{n\in\NN}\p(n/2)X^{|\co|n}
(\sum_{n\in2\NN}(\tit_n+\p(n/4)(wX)^n)^2\\&
=4\ti\La(wX)^2\prod_{\co\in\O_{\a,\g}}\prod_{k\ge1}(1-X^{|\co|2k})\i
=4\ti\La(wX)^2\prod_{k\ge1}\prod_{\co\in\O_{\a,\g}}(1-(X^{2k})^{|\co|})\i.\endalign$$
The product over $\co$ can be evaluated using 1.3. Using this and 3.16(e) we obtain
$$\align&\fc=4\ti\La(wX)^2\prod_{k\ge1}(1-qX^{4k})\i(1-X^{4k})^2
=4\ti\La(wX)^2\Ps_q(X^4)\Ps(X^4)^{-2}\\&
=\Ps_q(X^4)\Ps((wX)^2)^2(\Ps(-wX)\i+\Ps(wX)\i)^2\\&
=\Ps_q(X^4)\Ps((wX)^2)^2(\Ps(-wX)^{-2}+\Ps(wX)^{-2}+2\Ps(-wX)\i\Ps(X)\i).\endalign$$
Here we can replace $\Ps(wX)^{-2}+\Ps(-wX)^{-2}$ by 
$2\Ps((wX)^2)\i\Ps(-X^4)^{-2}\Ps(X^8)$.
(Indeed from $\sum_{j\in\ZZ}(-X)^{j^2}+\sum_{j\in\ZZ}X^{j^2}=2\sum_{j\in\ZZ}X^{4j^2}$ and
the Jacobi identity 13.6(b) we obtain
$$\Ps(X)^{-2}\Ps(X^2)+\Ps(-X)^{-2}\Ps(X^2)=2\Ps(-X^4)^{-2}\Ps(X^8).$$
We then replace $X$ by $wX$.) Thus we have
$$\align&\fc=2\Ps_q(X^4)\Ps((wX)^2)\Ps(-X^4)^{-2}\Ps(X^8)\\&
+2\Ps_q(X^4)\Ps((wX)^2)^2\Ps(-wX)\i\Ps(wX)\i.\endalign$$
The product of the last two factors above can be rewritten using 3.16(f) with $X$ replaced
by $wX$. We obtain
$$\fc=2\Ps_q(X^4)\Ps((wX)^2)\Ps(-X^4)^{-2}\Ps(X^8)+2\Ps_q(X^4)\Ps((wX)^2)\i\Ps(X^4).$$
Combining the results above we obtain
$$\align\sum_{n\in2\NN}&\a_nX^n=\fa+\fb'+\fb''+\fc=
3\Ps_q(X^4)\Ps((wX)^2)\Ps(-X^4)^{-2}\Ps(X^8)\\&+3\Ps_q(X^4)\Ps((wX)^2)\i\Ps(X^4)
+2\Ps_q(X^4)\Ps(X^2)\i\Ps(X^4).\tag c\endalign$$

\head 4. Counting data in the dual group\endhead
\subhead 4.1\endsubhead
Let $s\in\Xi$ be such that $\D:=\D_s\in\fA^N_0$. We show:

(a) ${}^\b\ps_s=\ps_{-s}\in\cs_{{}^\b\D}$.
\nl
It is enough to show that if $U\in\cu^\bul$ then
$[[\op_{\l\in U}V^{-s}_\l]]=N/2+[[\op_{\l\in U}V^s_\l]]$. This follows from 2.1(a) since
$\op_{\l\in U}V^{-s}_\l=\op_{\l\in-U}V^s_\l=\op_{\l\in U\i}V^s_\l$ is a subspace of $V$ 
complementary to $\op_{\l\in U}V^s_\l$.

\subhead 4.2\endsubhead
Let $\ufa$ (resp. $\ufa^N$) be the set of unordered pairs $(\D,\D')$ where $\D,\D'$ are 
elements of $\fA$ (resp. $\fA^N$) such that $\D'={}^\b\D$. For $k\in\{0,1,2\}$ let 
$\ufa^N_k$ be the set of all $(\D,\D')\in\ufa^N$ such that exactly $k$ of the numbers 
$n_1^\D,n_{-1}^\D$ are nonzero.

Let $\ufa^{N,\ne}$ (resp. $\ufa^{N,=}$) be the set of all $(\D,\D')\in\ufa^N$ such that
$\D\ne\D'$ (resp. $\D=\D'$). 

For $k\in\{0,1,2\}$ let $\ufa^{N,\ne}_k=\ufa^{N,\ne}\cap\ufa^N_k$, 
$\ufa^{N,=}_k=\ufa^{N,=}\cap\ufa^N_k$. Note that $\ufa^{N,=}_1=\em$.

Let $\tufa^N_0$ be the set of quadruples $\ti\he=(\D,\ps,\D',\ps')$ (with 
$(\D,\ps,\D',\ps')$, \lb $(\D',\ps',\D,\ps)$ identified) where $\D,\D'$ are elements of 
$\fA^N_0$ such that $\D'={}^\b\D$ and $\ps\in\cs_\D,\ps'\in\cs_{\D'}$ satisfy 
$\ps'={}^\b\ps$. (See 4.1(a).)

Let $\tufa^{N,\ne}_0$ (resp. $\tufa^{N,=}_0$) be the set of all
$(\D,\ps,\D',\ps')\in\tufa^N_0$ such that $\D\ne\D'$ (resp. $\D=\D'$). 

For any $\he\in\ufa^N$ we set $R_\he$ be the set of all $h\in\uXi$ such that
$\he=(\D_s,\D_{s'})$ where $\th\i(h)=\{s,s'\}$ (we have $s+s'=0$). For any 
$\ti\he\in\ufa^N_0$ we set $R_{\ti\he}$ be the set of all $h\in\uXi$ such that
$\ti\he=(\D_s,\ps_s,\D_{s'},\ps_{s'})$ where $\th\i(h)=\{s,s'\}$. Note that if 
$\D\in\ufa^N$ and $\he=(\D,{}^\b\D)$ then $R_\he=\th(R_\D)=\th(R_{{}^\b\D})$; if 
$\D\in\ufa^N_0,\ps\in\cs_\D$ and $\ti\he=(\D,\ps,{}^\b\D,{}^\b\ps)$ then 
$R_{\ti\he}=\th(R_\D^\ps)=\th(R_{{}^\b\D}^{{}^\b\ps})$. The following result and its 
proof are closely connected to 4.1(a). 

(a) {\it Let $\D\in\fA^N_0$. Assume that ${}^\b\D=\D$. Let $C,C'$ be two $SO_Q$-conjugacy
classes contained in $R_\D$. If $N\in2+4\NN$ then $C'=-C$. If $N\in4\NN$ then 
$-C=C,-C'=C'$.}
\nl
Let $s\in C$; we can find $x\in O(Q)$ such that $xsx\i=-s$. We have $xCx\i=C'$. It is 
enough to show that $\det x=(-1)^{N/2}$. Let $U\in\cu^\bul$. Then 
$E:=\op_{\l\in U}V^s_\l\in\tcs$ and $x(E)=\op_{\l\in\b(U)}V^s_\l\in\tcs$. Since 
$E\cap x(E)=0$ we see from 2.1(a) that $[[x(E)]]=N/2+[[E]]$ hence $\det(x)=(-1)^{N/2}$, 
as required.

In the setup of (a) and setting $\cs_\D=\{\ps_1,\ps_2\}$ we see that 
$\th(R_\D^{\ps_1})=\th(R_\D^{\ps_2})$ if $N\in2+4\NN$ and 
$\th(R_\D^{\ps_1})\ne\th(R_\D^{\ps_2})$ if $N\in4\NN$.

Using this and 2.2(b) we see that:

(b) {\it $R_\he$ ($\he\in\ufa^N_1\cup\ufa^N_2$), $R_\he$ ($\he\in\ufa^{N,=}_0$, 
$N\in2+4\NN$), $R_{\ti\he}$ ($\ti\he\in\tufa^{N,\ne}_0$, $R_{\ti\he}$ 
($\ti\he\in\tufa^{N,=}_0$, $N\in4\NN$) are exactly the $PSO_Q$-conjugacy classes in 
$\uXi$.}
\nl
For $h\in PSO_Q$ we denote by $Z(h)$ the centralizer of $h$ in $PSO_Q$. The component 
group $\ov{Z(h)}$ is isomorphic to:

$\ZZ/2\T\ZZ/2$ if $h\in R_\he,\he\in\ufa^{N,=}_2$, $N\in4\NN$;

$\ZZ/4$ if $h\in R_\he,\he\in\ufa^{N,=}_2$, $N\in2+4\NN$;

$\ZZ/2$ if $h\in R_\he,\he\in\ufa^{N,=}_0$, $N\in4\NN$;

$\ZZ/2$ if $h\in R_\he,\he\in\ufa^{N,\ne}_2$;

$\{1\}$ if $h\in R_\he,\he\in\ufa^{N,=}_0$, $N\in2+4\NN$;

$\{1\}$ if $h\in R_\he$, $\he\in\ufa^{N,\ne}_0\cup\ufa^{N,\ne}_1$.

\subhead 4.3\endsubhead
We define a permutation $\he\m{}^\g\he$ of $\ufa^N$ by $(\D,\D')\m({}^\g\D,{}^\g\D')$. The
fixed point set of this permutation is denoted by ${}^\g\ufa^N$. Note that for 
$(\D,\D')\in\ufa^N$ we have $(\D,\D')\in{}^\g\ufa^N$ if and only if ${}^\g\D$ is equal to 
$\D$ or to $\D'$ or equivalently if ${}^{\g_e}\D=\D$ for some $e\in\ZZ/2$.

Let $F\in\cf$. Clearly, $F:PSO_Q@>>>PSO_Q$ permutes the subsets $R_\he$ and a subset 
$R_\he$ is $F$-stable if and only if $\he\in{}^\g\ufa^N$. It follows that
$PSO_Q^F=\sqc_{\he\in{}^\g\ufa^N}R_\he^F$. We show:

(a) {\it Assume that $\ti\he=(\D,\ps,\D',\ps')\in\tufa^N_0$ satisfies 
${}^\g(\D,\D')=(\D,\D')$ so that ${}^{\g_e}\D=\D$ for some $e\in\ZZ/2$. Let 
$a=[V]_F+{}^ej_\D\in\ZZ/2$. Let $a+\ti\he=(\D,a+\ps,\D',a+\ps')$. Then 
$F(R_{\ti\he})=R_{a+\ti\he}$.}
\nl
Let $h\in R_{\ti\he}$. Let $s\in\th\i(h)$. We may assume that $\D=\D_s,\ps_s=\ps$. We have
$\D_{(-1)^eF(s)}={}^{\g_e}\D$ hence $\D_{(-1)^eF(s)}=\D$. By 2.5(a), 1.9(a), we have 
$\ps_{(-1)^eF(s)}(U)=a+\ps_s(U)$ for any $U\in\cu^\bst_e$. Hence 
$\ps_{(-1)^eF(s)}(U)=a+\ps_s(U)$ so that $\ps_{(-1)^eF(s)}=a+\ps$. We see that 
$F(h)=\th((-1)^eF(s))\in R_{a+\ti\he}$ and (a) follows. (Compare 2.7(a).)

Note that $e$ is not uniquely defined when $\D=\D'$. In that case and if $N\in4\NN$ we have
${}^0j_\D={}^1j_\D=0$, see 1.4(h), hence $a=[V]_F$ is independent of $e$; on the other hand
if $N\in2+4\NN$ then $R_{\ti\he}=R_{a+\ti\he}=R_\he$ where $\he=(\D,\D')$ hence (a) just
states that $F(R_\he)=R_\he$.

\subhead 4.4\endsubhead
For any $n\in\NN$ we fix $F:PGL_n(\kk)@>>>PGL_n(\kk)$, a Frobenius map for an
$\FF_q$-split rational structure and a diagram involution $\io:PGL_n@>>>PGL_n$ such that
$\io F=F\io$ and such that $\io\ne1$ whenever $n\ge3$. For any $\D\in\fA^N$ we set 
$$X^\D=\{(g_\l)\in\prod_{\l\in\kk^\tri}PGL_{n^\D_\l};g_\l=\io(g_{\l\i})\text{ for all }
\l\in\kk^\tri\}.$$
Let $e\in\ZZ/2$, $\D\in{}^{\g_e}\fA^N$. We define an automorphism $\mu_e:X^\D@>>>X^\D$ by 
$(g_\l)\m(g'_\l)$ where $g'_\l=g_{\g_e(\l)}$. This is well defined since $\g_e\a=\a\g_e$. 

Define $\fF:X^\D@>>>X^\D$ by $(g_\l)\m(g'_\l)$ where $g'_\l=F(g_\l)$; $\fF$ is a 
Frobenius map for an $\FF_q$-split rational structure on $X^\D$ commuting with $\mu_e$.

In the case where $\D\in\fA^N$ and ${}^\b\D=\D={}^\g\D$ we define an involution 
$\r:X^\D@>>>X^\D$ by $(g_\l)\m(\tg_\l)$ where $\tg_\l=g_{\b(\l)}$; for $e\in\ZZ/2$ we have 
$\r\mu_e=\mu_e\r$.

For any $n\in2\NN$ we write $Y_n$ instead of $(SO_{Q_n})_{\ad}$ where $Q_n$ is the quadratic
form $x_1x_2+x_3x_4+\do+x_{n-1}x_n$ on $\kk^n$.
We fix $\fF:Y_n@>>>Y_n$, a Frobenius map for an
$\FF_q$-split rational structure and a diagram involution 
$\nu:Y_n@>>>Y_n$ such that $\nu\fF=\fF\nu$ and such that $\nu\ne1$ whenever $n\ge4$. 
We define $\s:Y_n\T Y_n@>>>Y_n\T Y_n$ by $(g,g')\m(g',g)$.
For any $\D\in\fA^N$ we define $G^\D=X^\D\T Y_{n^\D_1}\T Y_{n^\D_{-1}}$.

\subhead 4.5\endsubhead
Let $e\in\ZZ/2$. For any semisimple conjugacy class $\cc$ in
$PSO_Q^{F_e}$ we choose $h\in\cc$ and $s\in\th\i(h)$. Let $\D=\D_s$, $\he=(\D,\D')$. Note 
that $\he\in{}^\g\ufa^N$. We can identify $Z(h)_{\ad}^0=G^\D$ (see 4.4) in such a way that 
$F_e:Z(h)_{\ad}^0@>>>Z(h)_{\ad}^0$ corresponds to $\fF\ph:G^\D@>>>G^\D$ where
$\fF:G^\D@>>>G^\D$ stands for $\fF\T\fF\T\fF$ of 4.4 (a Frobenius map for an $\FF_q$-split
rational structure on $G^\D$) and $\ph:G^\D@>>>G^\D$ is a diagram automorphism of $G^\D$ 
commuting with $\fF$. Let $\AA=\ov{Z(h)_{\ad}}^{F_e}$, a finite abelian $2$-group. Now any
element $f\in\AA$ gives rise to an automorphism $[f]$ of $Z(h)_{\ad}^0=G^\D$ (conjugation
by a representative of $f$ in $Z(h)_{\ad}^{F_e}$) which commutes with $\fF\ph$ and is well
defined up to composition with conjugation by an element in $(G^\D)^{\fF\ph}$. We now 
describe in each case the action 
of $\ph$ on $G^\D$. We also describe for each $f\in\AA$ a corresponding automorphism 
$[f]$ of $G^\D$ (it is enough to do that for a set of generators of $\AA$).

(I) Assume that $\he=(\D,\D')\in\ufa^{N,\ne}_0$. Define $e'\in\ZZ/2$ by 
${}^{\g_{e'}}\D=\D$. If $e\ne{}^{e'}j_\D$ then $R_\he^{F_e}=\em$. If $e={}^{e'}j_\D$ then 
$R_\he^{F_e}$ is a union of two $PSO_Q^{F_e}$-conjugacy classes $R_{\ti\he}^{F_e}$ where 
$\ti\he=(\D,\ps,\D',\ps')\in\tufa^N$. We have $\ph=\mu_{e'}\T1\T1$ and $\AA=\{1\}$.

(II) Assume that $\he=(\D,\D')\in\ufa^{N,\ne}_1$. We have ${}^\g\D=\D$ and $R_\he^{F_e}$ is
a single $PSO_Q^{F_e}$-conjugacy class. We have $\ph=\mu_0\T1\T\nu^{j_\D+e}$ (if 
$n^\D_{-1}\ne0$), $\ph=\mu_0\T\nu^{j_\D+e}\T1$ (if $n^\D_1\ne0$) and $\AA=\{1\}$.

(III) Assume that $\he=(\D,\D')\in\ufa^{N,\ne}_2$. Define $e'\in\ZZ/2$ by 
${}^{\g_{e'}}\D=\D$. Then $R_\he^{F_e}$ is a union of two $PSO_Q^{F_e}$-conjugacy classes 
$\cc_a (a\in\ZZ/2)$. For $\cc_a$ we have 
$\ph=\mu_{e'}\T\s^{e'}(\nu^a\T\nu^{a+e+{}^{e'}j_\D})$ and $\AA=\ZZ/2$ with generator $f$ 
such that $[f]=1\T\nu\T\nu$.

(IV) Assume that $\he=(\D,\D')\in\ufa^{N,=}_0$ and that $N\in2+4\NN$. We have 
${}^{\g_{e'}}\D=\D$ where $e'\in\ZZ/2$ is defined by $(q-(-1)^{e'})/2=e\in\ZZ/2$ and 
$R_\he^{F_e}$ is a single $PSO_Q^{F_e}$-conjugacy class. We have $\ph=\mu_{e'}\T1\T1$ and 
$\AA=\{1\}$. 

(V) Assume that $\he=(\D,\D')\in\ufa^{N,=}_0$ and that $N\in4\NN$. If $e=1$ then 
$R_\he^{F_e}=\em$. If $e=0$ then $R_\he^{F_e}$ is a union of four $PSO_Q^{F_e}$-conjugacy 
classes $\cc_{\ti\he,e'}$ where $\ti\he=(\D,\ps,\D',\ps')\in\tufa^N$ and $e'\in\ZZ/2$. For
$\cc_{\ti\he,e'}$ we have $\ph=\mu_{e'}\T1\T1$ and $\AA=\ZZ/2$ with generator $f$
such that $[f]=\r\T1\T1$.

(VI) Assume that $\he=(\D,\D')\in\ufa^{N,=}_2$ and that $N\in2+4\NN$. If $e=(q-1)/2\mod2$ 
then $R_\he^{F_e}$ is a union of four $PSO_Q^{F_e}$-conjugacy classes 
$\cc_{a,e'} (a\in\ZZ/2,e'\in\ZZ/2)$. For $\cc_{a,e'}$ we have
$\ph=\mu_{e'}\T\s^{e'}(\nu^a\T\nu^{\e'+a})$ and $\AA=\ZZ/4$ with generator $f$ such that
$[f]=\r\T\s(\nu\T1)$. If $e=(q+1)/2\mod2$ then $R_\he^{F_e}$ is a union of two 
$PSO_Q^{F_e}$-conjugacy classes $\cc_{e'} (e'\in\ZZ/2)$. For $\cc_{e'}$ we have 
$\ph=\mu_{e'}\T\s^{e'}(1\T\nu^{e'+1})$ and $\AA=\ZZ/2$ with generator $f$ such that 
$[f]=1\T\nu\T\nu$.

(VII) Assume that $\he=(\D,\D')\in\ufa^{N,=}_2$ and that $N\in4\NN$. If $e=0$ then 
$R_\he^{F_e}$ is a union of four $PSO_Q^{F_e}$-conjugacy classes 
$\cc_{a,e'} (a\in\ZZ/2,e'\in\ZZ/2)$. For $\cc_{a,e'}$ we have
$\ph=\mu_{e'}\T\s^{e'}(\nu^a\T\nu^a)$ and $\AA=\ZZ/2\T\ZZ/2$ with generators $f,f'$ such
that $[f]=\r\T\s$, $[f']=1\T\nu\T\nu$. If $e=1$ then $R_\he^{F_e}$ is a union of two 
$PSO_Q^{F_e}$-conjugacy classes $\cc_{e'} (e'\in\ZZ/2)$. For $\cc_{e'}$ we have 
$\ph=\mu_{e'}\T\s^{e'}(1\T\nu)$ and $\AA=\ZZ/2$ with generator $f$ such that
$[f]=1\T\nu\T\nu$.

\subhead 4.6\endsubhead
Let $\AA$ be a finite group acting on a finite set $\fS$. Let $\Irr_\AA(\fS)$ be the set of
isomorphism classes of irreducible $\AA$-equivariant local systems on $\fS$. 
For any integer $n\ge1$ let 

(a) $z_n=|\{x\in\fS;\text{ stabilizer of $x$ in $\AA$ has cardinal }n\}$.
\nl
Note that if $\AA$ is commutative then

(b) $|\Irr_\AA(\fS)|=|\AA|\i\sum_{n\ge1}n^2z_n.$
\nl
As an example, let $G$ be a reductive group with a Frobenius map $F:G@>>>G$ relative to an
$\FF_q$-rational structure. Let $\fU_F(G^0)$ be the set of {\it unipotent} representations
of $G^{0F}$ (up to isomorphism). Now any automorphism $\boa:G^0@>>>G^0$ commuting with $F$
induces a permutation of $\fU_F(G^0)$ denoted by $\un{\boa}$. In particular for any 
$g\in G^F$, $\Ad(g):G^{0F}@>>>G^{0F}$ induces a permutation 
$\un{\Ad(g)}:\fU_F(G^0)@>>>\fU_F(G^0)$ depending only on the image of $g$ in $(G/G^0)^F$. 
This defines an action of the finite group $(G/G^0)^F=G^F/G^{0F}$ on the finite set 
$\fU_F(G^0)$. We set
$$\fU_F(G)=\Irr_{(G/G^0)^F}(\fU_F(G^0)).$$

\subhead 4.7\endsubhead
For any $n\in\NN$ let $P_n$ be as in 3.9 so that $|P_n|=\p(n)$. If $F,\io$ are as in 4.4,
then we have canonically $\fU_F(PGL_n(\kk))=\fU_{F\io}(PGL_n(\kk))=P_n$. Let $e\in\ZZ/2$,
$\D\in{}^{\g_e}\fA^N$. Let $\mu_e:X^\D@>>>X^\D$ by as in 4.4. Write $\fU_e(X^\D)$ instead 
of $\fU_{\fF\mu_e}(X^\D)$. From the definitions we have canonically
$\fU_e(X^\D)=\prod_{\co\in\O_{\a,\g_e}}P_{n^\D_\co}$ hence
$$|\fU_e(X^\D)|=\vp_{\D,e}\tag a$$
where 
$$\vp_{\D,e}=\prod_{\co\in\O_{\a,\g_e}}\p(n^\D_\co).\tag b$$ 
Now assume that $\D\in\fA^N$ satisfies $\D={}^\b\D={}^\g\D$. We set
$$\vp_\D=\prod_{\co\in\O_{\a,\b,\g}}\p(n^\D_\co).\tag c$$
Note that $\r:X^\D@>>>X^\D$ induces a permutation 
$\un\r:\fU_e(X^\D)@>>>\fU_e(X^\D)$; from the definitions we have 
$$|(\fU_e(X^\D))^{\un\r}|=|(\prod_{\co\in\O_{\a,\g_e}}P_{n^\D_\co})^{\un\r}|=
=|\prod_{\co\in\O_{\a,\b,\g}}P_{n^\D_\co}|=\vp_\D.\tag d$$
Hence with respect to the $\ZZ/2$-action on $\fU_e(X^\D)$ given by $\un\r$ we have
$$\align&|\Irr_{\ZZ/2}(\fU_e(X^\D))|\\&=2|(\fU_e(X^\D))^{\un\r}|+(1/2)
|\fU_e(X^\D)-(\fU_e(X^\D))^{\un\r}|=(3/2)\vp_\D+(1/2)\vp_{\D,e}.\tag e\endalign$$

\subhead 4.8\endsubhead
For $n\in2\NN$ and $a\in\ZZ/2$ we write $\fU_a(Y_n)$ instead of $\fU_{\fF\nu^a}(Y_n)$; we
set 
$$\et_n=\p(n/4).\tag a$$ 
For $m\in\NN+1$ let $\Ph^0_m$ (resp. $\Ph^1_m$) be the set of "symbols" denoted by 
$\Ph_m^+$ (resp. $\Ph_m^-$) in \cite{\IRC, 3.3}. For $n\in2\NN+2,a\in\ZZ/2$ we can identify
$\fU_a(Y_n)=\Ph_{n/2}^a$, see \cite{\IRC}. For $n\in2\NN$ we define 
$$\x_n=|\Ph_{n/2}^1|\text{ if }n\ge2, \qua \x_0=0.\tag b$$
Then $|\Ph_{n/2}^0|=\x_n+2\et_n$ for $n\ge2$. Thus for $a\in\ZZ/2,n\ge2$ we have
$$|\Ph_{n/2}^a|=\x_n+(1+(-1)^a)\et_n.\tag c$$
The involution $\un{\nu}:\fU_a(Y_n)@>>>\fU_a(Y_n)$ induced by $\nu:Y_n@>>>Y_n$ becomes an
involution of $\Ph_{n/2}^a$ which is again denoted by $\un\nu$ (it can be also defined in 
terms of symbols). Note that $\un\nu$ acts as $1$ on $\Ph_{n/2}^1$ and
$|(\Ph_{n/2}^0)^{\un\nu}|=\x_n$. Hence with respect with the $\ZZ/2$ action on 
$\Ph_{n/2}^a$ given by $\un\nu$ we have
$$|\Irr_{\ZZ/2}(\Ph_{n/2}^a)|=2\x_n+(1/2)(1+(-1)^a)\et_n.\tag d$$
Similarly if $n,n'$ are in $2\NN+2$ then with respect to the $\ZZ/2$-action on
$\Ph_{n/2}^a\T\Ph_{n'/2}^{a'}$ given by $\nu\T\nu$ we have
$$\align&|\Irr_{\ZZ/2}(\Ph_{n/2}^a\T\Ph_{n'/2}^{a'})|=
(1/2)(4\x_n\x_{n'}+(1+(-1)^{a'})\x_n\et_{n'}\\&+
(1+(-1)^a)\et_n\x_{n'}+(1+(-1)^a)(1+(-1)^{a'})\et_n\et_{n'}).\tag e\endalign$$
For $n\in2\NN+2$ and $a,a'$ in $\ZZ/2$ we write $\fU_{a,a'}(Y_n\T Y_n)$ instead of \lb
$\fU_{(\fF\T\fF)\s(\nu^a\T\nu^{a'})}(Y_n\T Y_n)$. 
Then we have canonically
$$\fU_{a,a'}(Y_n\T Y_n)=\Ps_{n/2}^{a+a'}.\tag f$$

\subhead 4.9\endsubhead
Let $e\in\ZZ/2$. For any $\he\in{}^\g\ufa^N$ we set 
$$H^e_\he=\sum_\cc m_\cc$$
where $\cc$ runs over the set of $PSO_Q^{F_e}$-conjugacy classes contained in $R_\he^{F_e}$
and $m_\cc=|\fU_{F_e}(Z(h))|=|\fU_{F_e}(Z(h)_{\ad})|$ for some/any $h\in\cc$. We set
$$H_\he=H^0_\he-H^1_\he.$$
In 4.10-4.15 we compute $H_\he$ in the various cases (I)-(VII) in 4.5.

\subhead 4.10\endsubhead
In the setup of 4.5(I) we have $H^e_\he=2\vp_{\D,e'}$ if $e={}^{e'}j_\D$ and 
$H^e_\he=0$ if $e\ne{}^{e'}j_\D$. Hence
$$H_\he=(-1)^{{}^{e'}j_\D}2\vp_{\D,e'}.$$
In the setup of 4.5(II) we set $n=n^\D_1+n^\D_{-1}$ and we have
$$H^e_\he=\vp_{\D,0}|\Ph_{n/2}^{j_\D+e}|=\vp_{\D,0}(\x_n+(1+(-1)^{j_\D+e})\et_n).$$
hence
$$H_\he=(-1)^{j_\D}2\vp_{\D,0}\et_n.$$

\subhead 4.11\endsubhead
In the setup of 4.5(III) we set $n=n_1^\D,n'=n^\D_{-1}$,
$\x=\x_n,\x'=\x_{n'},\et=\et_n,\et'=\et_{n'},{}^{e'}j={}^0j_\D$. Assume first that $e'=0$.
Then 
$$\align&H_\he^e=\vp_{\D,0}\sum_{a,a';a+a'=e+j_\D}
|\Irr_{\ZZ/2}(\Ph_{n/2}^a\T\Ph_{n'/2}^{a'})|\\&=
(1/2)\vp_{\D,0}\sum\Sb a,a';\\a+a'=e+j_\D\endSb(4\x\x'+(1+(-1)^{a'})\x\et'\\&+
(1+(-1)^a)\et\x'+(1+(-1)^a)(1+(-1)^{a'})\et\et')\\&
=\vp_{\D,0}(4\x\x'+\x\et'+\et\x'+\et\et'+(-1)^{e+j_\D}\et\et'),\endalign$$
$$H_\he=(-1)^{j_\D}2\vp_{\D,0}\et\et'.$$
Assume next that $e'=1$. Then $n=n'$ and
$$H_\he^e=\vp_{\D,1}\sum_{a,a';a+a'=e+{}^1j_\D}|\Irr_{\ZZ/2}(\Ph_{n/2}^{a+a'})|=
2\vp_{\D,1}(2\x+(1/2)(1+(-1)^{e+{}^1j_\D})\et),$$
$$H_\he=(-1)^{{}^1j_\D}2\vp_{\D,1}\et.$$

\subhead 4.12\endsubhead
In the setup of 4.5(IV) we have $H^e_\he=\vp_{\D,e'}$ where $e'$ is defined by 
$(q-(-1)^{e'})/2=e\in\ZZ/2$. Hence
$$H_\he=(-1)^{(q-1)/2}(\vp_{\D,0}-\vp_{\D,1}).$$

\subhead 4.13\endsubhead
In the setup of 4.5(V) we have (using 4.7(e)):
$$H_\he^0=(3\vp_\D+\vp_{\D,0})+(3\vp_\D+\vp_{\D,1})$$
and $H_\he^1=0$. Hence
$$H_\he=6\vp_\D+\vp_{\D,0}+\vp_{\D,1}.$$

\subhead 4.14\endsubhead
In the setup of 4.5(VI) we set $n=n^\D_1=n^\D_{-1}$, $\x=\x_n,\et=\et_n$. We have
$$H_\he^e=\sum_{a\in\ZZ/2}|\Irr_{\ZZ/4}(\fU_0(X^\D)\T\Ph_{n/2}^a\T\Ph_{n/2}^a)|
+2|\Irr_{\ZZ/4}(\fU_1(X^\D)\T\Ph_{n/2}^1)|$$
if $e=(q-1)/2\in\ZZ/2$,
$$H_\he^e=|\Irr_{\ZZ/2}(\fU_0(X^\D)\T\Ph_{n/2}^0\T\Ph_{n/2}^1)|+
|\Irr_{\ZZ/2}(\fU_1(X^\D)\T\Ph_{n/2}^0)|$$
if $e=(q+1)/2\in\ZZ/2$, where the actions of $\ZZ/4,\ZZ/2$ are as below. (For a 
$\ZZ/4$-action on a set $\fS$ we 
denote by $z_i$ the number of elements in $\fS$ whose stabilizer in $\ZZ/4$ has cardinal 
$i$.) The $\ZZ/4$ action on $\fU_0(X^\D)\T\Ph_{n/2}^0\T\Ph_{n/2}^0$ is such that a 
generator acts by $(A,B,C)\m(\un\r(A),C,\nu(B))$. Then $(A,B,C)$ is fixed by the generator
if $\un\r A=A$, $B=\un\nu B=C$; it is fixed by the square of the generator if 
$B=\un\nu B$, $C=\un\nu C$. We see that
$$z_4=\vp_\D\x,\qua z_2=-\vp_\D\x+\vp_{\D,0}\x^2,\qua z_1=4\vp_{\D,0}(\et^2+\x\et).$$
We have 
$$|\Irr_{\ZZ/4}(\fU_0(X^\D)\T\Ph_{n/2}^0\T\Ph_{n/2}^0)|=
4z_4+z_2+(1/4)z_1=3\vp_\D\x+\vp_{\D,0}(\x^2+\et^2+\x\et).$$
The same computation shows that
$$|\Irr_{\ZZ/4}(\fU_0(X^\D)\T\Ph_{n/2}^1\T\Ph_{n/2}^1)|=
3\vp_\D\x+\vp_{\D,0}\x^2.$$
(We put formally $\et=0$ in the previous computation.)

The $\ZZ/4$ action on $\fU_1(X^\D)\T\Ph_{n/2}^1)$ is such that a generator acts by 
$(A,B)\m(\un\r(A),\nu(B))$. Then $(A,B)$ is fixed by the generator if $\un\r A=A$; it is 
always fixed by the square of the generator. We see that 
$$z_4=\vp_\D\x,\qua z_2=\vp_{\D,1}\x-\vp_\D\x,\qua z_1=0.$$
We have 
$$|\Irr_{\ZZ/4}(\fU_1(X^\D)\T\Ph_{n/2}^1)|=4z_4+z_2+(1/4)z_1=3\vp_\D\x+\vp_{\D,1}\x.$$
We see that for $e=(q-1)/2\mod2$ we have
$$\align&H_\he^e=3\vp_\D\x+\vp_{\D,0}(\x^2+\et^2+\x\et)+3\vp_\D\x+\vp_{\D,0}\x^2
+2(3\vp_\D\x+\vp_{\D,1}\x)\\&=12\vp_\D\x+\vp_{\D,0}(2\x^2+\et^2+\x\et)+2\vp_{\D,1}\x.
\endalign$$
The $\ZZ/2$ action on $\fU_0(X^\D)\T\Ph_{n/2}^0\T\Ph_{n/2}^1$ is by
$(A,B,C)\m(A,\nu(B),\nu(C))$. Hence
$$\align&|\Irr_{\ZZ/2}(\fU_0(X^\D)\T\Ph_{n/2}^0\T\Ph_{n/2}^1)|=
\vp_{\D,0}|\Irr_{\ZZ/2}(\Ph_{n/2}^0\T\Ph_{n/2}^1)|\\&=\vp_{\D,0}(2\x^2+\x\et).\endalign$$
(The last equality comes from 4.8(e).) The $\ZZ/2$-action on $\fU_1(X^\D)\T\Ph_{n/2}^0$ is
by $(A,B)\m(A,\nu(B))$. Hence
$$|\Irr_{\ZZ/2}(\fU_1(X^\D)\T\Ph_{n/2}^0)|=\vp_{\D,1}(2\x+\et).$$
We see that for $e=(q+1)/2\mod2$ we have
$$H^e_\he=\vp_{\D,0}(2\x^2+\x\et)+\vp_{\D,1}(2\x+\et).$$
It follows that
$$H_\he=(-1)^{(q-1)/2}(12\vp_\D\x+\vp_{\D,0}\et^2-\vp_{\D,1}\et).$$

\subhead 4.15\endsubhead
In the setup of 4.5(VII) we set $n=n^\D_1=n^\D_{-1}$, $\x=\x_n,\et=\et_n$. We have
$$H_\he^0=\sum_{a\in\ZZ/2}|\Irr_{\ZZ/2\T\ZZ/2}(\fU_0(X^\D)\T\Ph_{n/2}^a\T\Ph_{n/2}^a)|
+2|\Irr_{\ZZ/2\T\ZZ/2}(\fU_1(X^\D)\T\Ph_{n/2}^0)|,$$
$$H_\he^1=|\Irr_{\ZZ/2}(\fU_0(X^\D)\T\Ph_{n/2}^0\T\Ph_{n/2}^1)|+
|\Irr_{\ZZ/2}(\fU_1(X^\D)\T\Ph_{n/2}^1)|,$$
where the actions of $\ZZ/2\T\ZZ/2,\ZZ/2$ are as below. (For a $\ZZ/2\T\ZZ/2$-action on a 
set $\fS$ we denote by $z_i$ the number of elements in $\fS$ whose stabilizer in 
$\ZZ/2\T\ZZ/2$ has cardinal $i$.)
The $\ZZ/2\T\ZZ/2$ action on $\fU_0(X^\D)\T\Ph_{n/2}^0\T\Ph_{n/2}^0$ is such that some
element $f$ acts as $(A,B,C)\m(\un\r(A),C,B)$ and some element $f'$ acts as
$(A,B,C)\m(A,\un\nu(B),\un\nu(C))$.
An element $(A,B,C)$ is fixed by $f$ if $A=\un\r(A),B=C$; it is fixed by $f'$ if
$B=\un\nu(B),C=\un\nu(C)$; it is fixed by $ff'$ if $A=\un\r(A),B=\un\nu(C)$. We see that
$$z_4=\vp_\D\x,\qua z_2=\vp_\D(-\x+4\et)+\vp_{\D,0}\x^2,\qua
z_1=4\vp_{\D,0}(\et^2+\x\et)-4\vp_\D\et.$$
We have
$$\align&|\Irr_{\ZZ/2\T\ZZ/2}(\fU_0(X^\D)\T\Ph_{n/2}^0\T\Ph_{n/2}^0)|\\&
=4z_4+z_2+(1/4)z_1=3\vp_\D(\x+\et)+\vp_{\D,0}(\x^2+\et^2+\x\et).\endalign$$
The same computation shows that
$$|\Irr_{\ZZ/2\T\ZZ/2}(\fU_0(X^\D)\T\Ph_{n/2}^1\T\Ph_{n/2}^1)|=3\vp_\D\x+\vp_{\D,0}\x^2.$$
(We put formally $\et=0$ in the previous computation.)
The $\ZZ/2\T\ZZ/2$-action on $\fU_1(X^\D)\T\Ph_{n/2}^0$ is such that $f$ acts as 
$(A,B)\m(\un\r(A),B)$ and $f'$ acts as $(A,B)\m(A,\un\nu(B))$.
An element $(A,B)$ is fixed by $f$ if $A=\un\r(A)$; it is fixed by $f'$ if
$B=\un\nu(B)$; it is fixed by $ff'$ if $A=\un\r(A),B=\un\nu(C)$; it is fixed by 
$\ZZ/2\T\ZZ/2$ if $A=\un\r(A),B=\un\nu(C)$. We see that
$$z_4=\vp_\D\x,\qua z_2=\vp_\D(-\x+2\et)+\vp_{\D,1}\x,\qua z_1=-2\vp_\D\et+2\vp_{\D,1}\et.
$$
Hence
$$\align&|\Irr_{\ZZ/2\T\ZZ/2}(\fU_1(X^\D)\T\Ph_{n/2}^0)|\\&=
4z_4+z_2+(1/4)z_1=(3/2)\vp_\D(2\x+\et)+(1/2)\vp_{\D,1}(2\x+\et).\endalign$$
We see that 
$$\align&H_\he^0=3\vp_\D(\x+\et)+\vp_{\D,0}(\x^2+\et^2+\x\et)\\&+
3\vp_\D\x+\vp_{\D,0}\x^2+3\vp_\D(2\x+\et)+\vp_{\D,1}(2\x+\et)\\&=
6\vp_\D(2\x+\et)+\vp_{\D,0}(2x^2+\et^2+\x\et)+\vp_{\D,1}(2\x+\et).\endalign$$
The $\ZZ/2$ action on $\fU_0(X^\D)\T\Ph_{n/2}^0\T\Ph_{n/2}^1$ is by
$(A,B,C)\m(A,\nu(B),\nu(C))$. As in 4.14 we have
$$|\Irr_{\ZZ/2}(\fU_0(X^\D)\T\Ph_{n/2}^0\T\Ph_{n/2}^1)|=
\vp_{\D,0}(2\x^2+\x\et).$$
The $\ZZ/2$ action on $\fU_1(X^\D)\T\Ph_{n/2}^1$ is trivial. Hence
$$|\Irr_{\ZZ/2}(\fU_1(X^\D)\T\Ph_{n/2}^1)|=2\vp_{\D,1}\x.$$
We see that
$$H_\he^1=\vp_{\D,0}(2\x^2+\x\et)+2\vp_{\D,1}\x.$$
It follows that
$$H_\he=6\vp_\D(2\x+\et)+\vp_{\D,0}\et^2+\vp_{\D,1}\et.$$

\subhead 4.16\endsubhead
For any $\D\in{}^\g\fA^N$ we set $\hat{\vp}_{\D,0}=\vp_{\D,0}\et_{n^\D_1}\et_{n^\D_{-1}}$.
For any $\D\in{}^{\g_1}\fA^N$ we set $\hat{\vp}_{\D,1}=\vp_{\D,1}\et_{n^\D_1}$. For any 
$\D\in\fA^N$ we set $u_\D=1$ if ${}^\b\D=\D$ and $u_\D=2$ if ${}^\b\D\ne\D$. For any 
$\D\in\fA^N$ we set

$\hat{\vp}_\D=\vp_\D(2\x_{n^\D_1}+(1/2)(1+(-1)^{n^\D_J})\et_{n^\D_1})$ 
if $\D={}^\b\D={}^\g\D$,

$\hat{\vp}_\D=0$, otherwise.
\nl
(Recall that $J$ is as in 1.1.)

Let $\he=(\D,\D')\in{}^\g\ufa^N$. We have
$$H_\he=(-1)^{n^\D_J(q-1)/2}6\hat{\vp}_\D
+u_\D\sum_{e'\in\ZZ/2;{}^{\g_{e'}}\D=\D}(-1)^{{}^{e'}j_\D}\hat{\vp}_{\D,e'}.$$
This is verified by considering the various cases in 4.10-4.15; we also make use of 1.4(h)
and the equality $N/2=n^\D_J\in\ZZ/2$ from the proof of 1.4(h).

\subhead 4.17\endsubhead
Let $e\in\ZZ/2$. We set 
$$\ha^e_N=\sum_\cc m_\cc\tag a$$
where $\cc$ runs over the set of all semisimple $PSO_Q^{F_e}$-conjugacy classes in 
$PSO_Q^{F_e}$ and $m_\cc$ is as in 4.9. We set 
$$\ha_N=\ha^0_N-\ha^1_N.$$
From the definitions we have $\ha_N=\sum_{\he\in{}^\g\ufa^N}H_\he$. From this and 4.16 we see 
that
$$\sum_{n\in2\NN}\ha_nX^n=\fa_0+\fa_1+6\fd\tag b$$
where $\ha_0=8$,
$$\fa_{e'}=\sum_{\D\in{}^{\g_{e'}}\fA}(-1)^{{}^{e'}j_\D}\hat{\vp}_{\D,e'}X^{m_\D},$$
$$\fd=\sum_{\D\in\fA;{}^\b\D=\D={}^\g\D}(-1)^{n^\D_J(q-1)/2}\hat{\vp}_\D X^{m_\D}$$
and $m_\D=n$ for $\D\in\fA^n$.

\subhead 4.18\endsubhead
For $e\in\ZZ/2$ let $y_e=\sum_{j\in\ZZ}((-1)^eX^4)^{j^2}$. By Jacobi's identity we have

(a) $y_{-1}=\Ps(X^4)^{-2}\Ps(X^8),\qua y_1=\Ps(-X^4)^{-2}\Ps(X^8)$.
\nl
($\Ps$ as in 3.16.) By \cite{\IRC, 3.4.2)} we have 

(b) $\sum_{n\ge0}(\x_n+\et_n)X^{2n}=(1/4)\Ps(X^4)^2(y_1+3y_{-1})$,

(c) $\sum_{n\ge0}\et_nX^{2n}=\Ps(X^4)^2y_{-1}$.
\nl
It follows that

(d) $\sum_{n\ge0}(2\x_n+(1/2)\et_n)X^{2n}=(1/2)\Ps(X^4)^2y_1$.

\subhead 4.19\endsubhead
We set 
$$\fd'=\sum_{n\in2\NN,n'\in\NN}\uu^{n'}(2\x_n+(1/2)(1+(-1)^{n'})\et_n)
\p(n')X^{2n+2n'}.$$
We have
$$\align&
\fd'=\sum_{n'\in\NN}\uu^{n'}\p(n')X^{2n'}\sum_{n\in2\NN}(2\x_n+(1/2)\et_n)X^{2n}\\&
+(1/2)\sum_{n'\in\NN}(-\uu)^{n'}\p(n')X^{2n'}\sum_{n\in2\NN}\et_nX^{2n}\\&
=\Ps(\uu X^2)\sum_n(2\x_n+(1/2)\et_n)X^{2n}
+(1/2)\Ps(-\uu X^2)\sum_n\et_nX^{2n}.\endalign$$
Using now 4.18(d),(c),(a) we obtain
$$\align&\fd'=(1/2)\Ps(\uu X^2)\Ps(X^4)^2y_1+(1/2)\Ps(-\uu X^2)\Ps(X^4)^2y_{-1}\\&
=(1/2)\Ps(\uu X^2)\Ps(X^4)^2\Ps(-X^4)^{-2}\Ps(X^8)+(1/2)\Ps(-\uu X^2)\Ps(X^8).\endalign$$

\subhead 4.20\endsubhead
We compute the sum $\fd$ in 4.17. We have
$$\align&\fd=\sum_{\D\in\fA;{}^\b\D=\D={}^\g\D}(-1)^{n^\D_J(q-1)/2}
\prod_{\co\in\O_{\a,\b,\g}}\p(n_\co^\D)X^{|\co|n^\D_\co}\\&
\T(2\x_{n^\D_1}+(1/2)(1+(-1)^{n_J})\et_{n^\D_1})X^{2n^D_1}
=\fd'\prod_{\co\in\O_{\a,\b,\g};\co\ne J}(\sum_{n\in\NN}\p(n)X^{|\co|n})\\&
=\fd'\prod_{\co\in\O_{\a,\b,\g};\co\ne J}\prod_{k\ge1}(1-X^{|\co|k})\i
=\fd'\prod_{k\ge1}\prod_{\co\in\O_{\a,\b,\g};\co\ne J}(1-X^{|\co|k})\i.\endalign$$
The product over $\co$ can be evaluated using 1.6(a). We obtain 
$$\fd=\fd'\prod_{k\ge1}\{(1-qX^{4k})\i(1-X^{4k})^2\}=\Ps_q(X^4)\Ps(X^4)^{-2}\fd'.$$
Introducing here the formula for $\fd'$ at the end of 4.19, we obtain
$$\fd=(\Ps_q(X^4)\Ps(\uu X^2)\Ps(-X^4)^{-2}\Ps(X^8)
+\Ps_q(X^4)\Ps(-\uu X^2)\Ps(X^4)^{-2}\Ps(X^8))/2.$$
Here we substitute $\Ps(-\uu X^2)\Ps(X^4)^{-2}\Ps(X^8)=\Ps(\uu X^2)\i\Ps(X^4)$ which 
follows  from 3.16(f) with $X$ replaced by $\uu X^2$. We obtain
$$\fd=(1/2)\Ps_q(X^4)\Ps(\uu X^2)\Ps(-X^4)^{-2}\Ps(X^8)
+(1/2)\Ps_q(X^4)\Ps(\uu X^2)\i\Ps(X^4).\tag a$$

\subhead 4.21\endsubhead
We have
$$\fa_0=\sum_{\D\in{}^\g\fA}(-1)^{j_\D}\prod_{\co\in\O_{\a,\g}}\p(n^\D_\co)
X^{|\co|\n^\D_\co}\et_{n^\D_1}X^{n^\D_1}\et_{n^\D_{-1}}X^{n^\D_{-1}},$$
$$\fa_1=\sum_{\D\in{}^{\g_1}\fA}(-1)^{{}^1j_\D}\prod_{\co\in\O_{\a,\g_1}}\p(n^\D_\co)
X^{|\co|\n^\D_\co}\et_{n^\D_1}X^{2n^\D_1}.$$
Thus
$$\fa_0=\prod_{\co\in\O_{\a,\g}}(\sum_{n\in\NN}\p(n)(-1)^{j(\co)n}X^{|\co|n})
(\sum_{n\in\NN}\p(n/4)X^n)^2,$$
$$\fa_1=\prod_{\co\in\O_{\a,\g_1}}(\sum_{n\in\NN}\p(n)(-1)^{{}^1j(\co)n}X^{|\co|n})
(\sum_{n\in\NN}\p(n/4)X^{2n}),$$
$$\align&\fa_0=
\prod_{\co\in\O_{\a,\g}}\prod_{k\ge1}(1-((-1)^{j(\co)k}X^{|\co|k})\i\Ps(X^4)^2\\&
=\prod_{k\ge1;k\text{ odd}}\prod_{\co\in\O_{\a,\g}}(1-(-1)^{j(\co)}X^{|\co|k})\i
\prod_{k\ge1;k\text{ even}}\prod_{\co\in\O_{\a,\g}}(1-X^{|\co|k})\i\Ps(X^4)^2,\endalign$$
$$\align&\fa_1
=\prod_{\co\in\O_{\a,\g_1}}\prod_{k\ge1}(1-((-1)^{{}^1j(\co)k}X^{|\co|k})\i\Ps(X^8)\\&
=\prod_{k\ge1;k\text{ odd}}\prod_{\co\in\O_{\a,\g_1}}(1-(-1)^{{}^1j(\co)}X^{|\co|k})\i
\prod_{k\ge1;k\text{ even}}\prod_{\co\in\O_{\a,\g_1}}(1-X^{|\co|k})\i\Ps(X^8).\endalign$$
In the last formulas for $\fa_0,\fa_1$ we replace the products over $\co$ by the
expressions provided by 1.4(e),(f) and 1.3 (with $X$ replaced by $X^k$). We obtain
$$\align&\fa_0=\prod_{k\ge1;k\text{ odd}}(1-X^{2k})\prod_{k\ge1;k\text{ even}}
\{(1-qX^{2k})\i(1-X^{2k})^2\}\Ps(X^4)^2\\&=\Ps(X^2)\i\Ps(X^4)\Ps_q(X^4),\endalign$$
$$\align&\fa_1=\prod_{k\ge1;k\text{ odd}}(1-X^{2k})\prod_{k\ge1;k\text{ even}}
\{(1-qX^{2k})\i(1-X^{4k})\}\Ps(X^8)\\&=\Ps(X^2)\i\Ps(X^4)\Ps_q(X^4).\endalign$$
Thus we have
$$\fa_0=\fa_1=\Ps(X^2)\i\Ps(X^4)\Ps_q(X^4).$$
Combining this with 4.17(b), 4.20(a) we see that
$$\align&\sum_{n\in2\NN}\ha_nX^n=3\Ps_q(X^4)\Ps(\uu X^2)\Ps(-X^4)^{-2}\Ps(X^8)\\&+
3\Ps_q(X^4)\Ps(\uu X^2)\i\Ps(X^4)+2\Ps_q(X^4)\Ps(X^2)\i\Ps(X^4).\endalign$$
Since this equality has the same right hand side as 3.17(c), we deduce:
$$\sum_{n\in2\NN}\ha_nX^n=\sum_{n\in2\NN}\a_nX^n.$$
It follows that $\ha_N=\a_N$ that is,
$$\ha_N^0-\ha_N^1=|\Irr(\tSO_Q^{F_0})|-|\Irr(\tSO_Q^{F_1})|.\tag a$$

\head 5. Classification of irreducible representations\endhead
\subhead 5.1\endsubhead
Let $G$ be a connected reductive group defined over $\FF_q$ with Frobenius map $F:G@>>>G$.
Let $G^*$ be a connected reductive group defined over $\FF_q$, dual to $G$ as in 
\cite{\DL}. We denote again by $F:G^*@>>>G^*$ the corresponding Frobenius map. For any 
$F$-stable maximal torus $T$ of $G^*$ and any $s\in T^F$ let $R_T^G(s)$ be the virtual 
representation of $G^F$ attached in \cite{\DL} to an $F$-stable maximal torus of $G$ 
corresponding to $T$ and to a character of the group of its rational points corresponding 
to $s$. 

Let $\fs(G^{*F})$ be the set of semisimple $G^{*F}$-conjugacy classes in $G^{*F}$. For 
$\cc\in\fs(G^{*F})$ let $\Irr^\cc(G^F)$ be the set of all $r\in\Irr(G^F)$ which appear 
with nonzero multipicity in $R_T^G(s)$ for some $s\in\cc$ and some $F$-stable maximal torus
$T$ of $G^*$ such that $s\in T$. We have a partition 
$$\Irr(G^F)=\sqc_{\cc\in\fs(G^{*F})}\Irr^\cc(G^F),\tag a$$
see \cite{\DL}, \cite{\IRC}. In particular, $\Irr^1(G^F)=\fU_F(G)$ (notation of 4.6).

For a semisimple element $s\in G^*F$ let $Z(s)$ be the centralizer of $s$ in $G^*$. Assume
that $\cc\in\fs(G^{*F})$. 
Now $G^{*F}$ acts in an obvious way on $\sqc_{s\in\cc}\fU_F(Z(s)^0)$ inducing the
conjugation action on $\cc$. Let $\fU_\cc=\Irr_{G^F}(\sqc_{s\in\cc}\fU_F(Z(s)^0))$
(see 4.6). Note that for any $s\in\cc$ there is an obvious restriction isomorphism
$\fU_\cc@>\si>>\fU_F(Z(s))$.

Let $\e_G$ be $1$ if the $\FF_q$-rank of $G$ is even and $\e_G=-1$, otherwise.
We can state the following result.

\proclaim{Theorem 5.2} Let $\cc$ be a semisimple $G^{*F}$-conjugacy class in $G^{*F}$.
There exists a bijection 
$$\Irr^\cc(G^F)\lra\fU_\cc\tag a$$ 
such that the following holds. Let $s\in\cc$ and let $\r\in\Irr_\cc(G^F)$, 
$\r'\in\fU_\cc=\fU_F(Z(s))$ correspond to each other under (a). Then $\r'$ determines a
$\ov{Z(s)}^F$-orbit $\{\r'_1,\r'_2,\do,\r'_k\}$ in $\fU_F(Z(s)^0)$. Let $T$ be a maximal 
torus of $Z(s)^0$ which is $F$-stable. Then the multiplicity of $\r$ in $R_T^G(s)$ is 
$\e_G\e_{Z(s)^0}$ times the sum over $i\in[1,k]$ of the multiplicities of $\r'_i$ in 
$R_T^{Z(s)^0}(1)$.
\endproclaim
In \cite{\OR} it is shown that the theorem holds when $G$ has connected centre. In
\cite{\AST} it is shown that the theorem holds without assumption on the centre of $G$ if
we assume that a certain multiplicity one statement holds for $\tSO_Q$ with $N\in4\NN$.
(The reason that $\tSO_Q$ is especially difficult is that it is the only almost simple
group whose center is noncyclic.) In this section we show how to prove the required 
multiplicity one statement.

\subhead 5.3\endsubhead
In the remainder of this section we assume that $N\in4+4\NN$. Recall that 
$\k\i\{1,-1\}$ is equal to
the centre of $\tSO_Q$; it is isomorphic to $\ZZ/2\T\ZZ/2$. Let $\ct$ be the product of two
copies of $\kk^*$ indexed by the elements of $\k\i(-1)$. There is a unique homomorphism
$g\m t_g$, $\k\i\{1,-1\}@>>>\ct$ such that for $\o\in\k\i(-1)$ the component of $t_\o$ in 
the copy of $\kk^*$ indexed by $\o$ (resp. $-\o$) has order $2$ (resp. is $1$). This
identifies $k\i\{1,-1\}$ with the subgroup $\{t\in\ct;t^2=1\}$. Let $G$ be the 
quotient of $\tSO_Q\T\ct$ by $\k\i\{1,-1\}$, imbedded diagonally in $\tSO_Q\T\ct$. We 
identify $\tSO_Q,\ct$ with closed subgroups of $G$ via the obvious homomorphisms
$\tSO_Q@>>>G$, $\ct@>>>G$. Note that $\tSO_Q$ is the derived subgroup of
$G$ and $\ct$ is the center of $G$; moreover, $G/\ct=PSO_Q$.

By 2.8(a) with $y=-1$, for $e\in\ZZ/2$ and $\o\in\k\i(-1)$ we have $F_e(\o)=(-1)^e\o$. We 
define $F_e:\ct@>>>\ct$ as $1$ if $e=0$ and as $(x,x')\m(x'{}^q,x^q)$ if $e=1$. Then the
imbedding $\k\i\{1,-1\}\sub\ct$ is compatible with the actions of $F_e$. Hence there is a
unique $\FF_q$-rational structure on $G$ such that the corresponding Frobenius map
$G@>>>G$ maps $gt$ to $F_e(g)F_e(t)$ for any $g\in\tSO_Q$ and any $t\in\ct$;
this Frobenius map is denoted again by $F_e$.

Let $C=G/\tSO_Q$ so that $C^{F_e}=G^{F_e}/\tSO_Q^{F_e}$.
We show that $C^{F_1}$ is a cyclic group. We can identify $C=\kk^*\T\kk^*$ so that
$F_1:C@>>>C$ becomes $(x,x')\m(x'{}^q,x^q)$. Then $C^{F_1}\cong\FF_{q^2}^*$ so that $C^{F_1}$ is
indeed cyclic. Hence any
irreducible representation of $G^{F_1}$ is multiplicity free when restricted to
$\tSO_Q^{F_1}$ (see \cite{\AST, 9(b)}). By the arguments in \cite{\AST, 11} we see that
Theorem 5.2 holds for $\tSO_Q^{F_1}$. It follows that $|\Irr(\tSO_Q^{F_1})|=\ha_N^1$ where 
$\ha_N^1$ is as in 4.17(a). Using this and 4.21 we deduce:
$$|\Irr(\tSO_Q^{F_0})|=\ha_N^0.\tag a$$

\subhead 5.4\endsubhead
Let $G$ be as in 5.3. In the remainder of this section we take $F=F_0$. We have naturally 
$PSO_Q=G^*/\ct'$ (as groups defined 
over $\FF_q$) where $\ct'$ is the centre of $G^*$ (a two-dimensional torus). Moreover, the
derived subgroup of $G^*$ is $\tSO_Q$. Lt $\p:G^*@>>>PSO_Q$ be the canonical map.

Since $G^*$ has connected centre, Theorem 5.2 holds for $G^F$. In
particular, for any $\cc\in\fs(G^{*F})$ there exists a bijection
$$\Irr^\cc(G^F)\lra\fU_\cc\tag a$$
with the following property: if $s\in\cc$, $\r\in\Irr_\cc(G^F)$,
$\r'\in\fU_\cc=\fU_F(Z(s))$ correspond to each other under (a) and $T$ is a maximal 
torus of $Z(s)$ which is $F$-stable, then the multiplicity of $\r$ in $R_T^G(s)$ is 
$\e_G\e_{Z(s)}$ times the multiplicity of $\r'$ in $R_T^{Z(s)}(1)$. Note that in our case 
the property above determines the bijection (a) uniquely.

If $t\in\ct'{}^F$ and $\cc\in\fs(G^{*F})$ then $t\cc\in\fs(G^{*F})$ and 
multiplication by $t$ induces an isomorphism $\fU_\cc@>>>\fU_{t\cc}$ (this comes
from an isomorphism $\fU_F(Z(s))@>>>\fU_F(Z(ts))$ induced by
$t:Z(s)@>>>Z(ts)$ for any $s\in\cc$). Thus we obtain an action of $\ct'{}^F$ on the set
$\sqc_{\cc\in\fs(G^{*F})}\fU_\cc$, hence via (a) and 5.1(a), an action of $\ct'{}^F$ on 
the set $\Irr(G^F)$. We have $\ct'{}^F=C_0^*$ where 
$C_0^*=\Hom(C^F,\bbq^*)$. Thus we obtain an action of the abelian group $C^*_0$ on 
$\Irr(G^F)$. From the definitions we see that this coincides with the action of $C^*_0$
given by tensoring by a one dimensional representation.

We show that the stabilizer in $C^*_0=\ct'{}^F$ of any $r\in\Irr(G^F)$ has order
dividing $4$. This stabilizer is contained in $\ct'_\cc:=\{t\in\ct'{}^F;t\cc=\cc\}$ for
some $\cc\in\fs(G^{*F})$. If $t\in\ct'_\cc$ then for $s\in\cc$ we have $ts=gsg\i$ for
some $g\in G^{*F}$; in particular $t$ is in the derived subgroup of $G^*$ hence it is in
the centre of $\tSO_Q$ so that $|\ct'_\cc|$ divides $4$.

For any $i$ let $\tz_i$ be the number of $r\in\Irr(G^F)$ whose stabilizer in
$C^*_0=\ct'{}^F$ has order $i$. Let $j:\fs(G^{*F})@>>>\fs(PSO_Q^F)$ be the map
induced by $\p$. For any $\cc'\in\fs(PSO_Q^F)$, the subset 
$\sqc_{\cc\in j\i(\cc')}\fU_\cc$ of $\sqc_{\cc\in\fs(G^{*F})}\fU_\cc$ is stable
under the $\ct'{}^F$-action; we denote by $\tz_{i,\cc'}$ the number of elements in
$\sqc_{\cc\in j\i(\cc')}\fU_\cc$ whose stabilizer in $\ct'{}^F$ has order $i$. 
Note that
$$\tz_i=\sum_{\cc'\in\fs(PSO_Q^F)}\tz_{i,\cc'}.\tag b$$
For $\cc'\in\fs(PSO_Q^F)$ we choose $\cc\in j\i(\cc')$, $s\in\cc$ and we let
$s'=\p(s)\in\cc'$. Since $\ct'{}^F$ acts transitively on $j\i(\cc')$ and the isotropy 
group of $\cc$ is $\ct'{}^F$, we have
$$\tz_{i,\cc'}=|\ct'{}^F||\ct'_\cc|\i\tz_{i,\cc}\tag c$$
where $\tz_{i,\cc}$ is the number of elements in $\fU_\cc$ whose stabilizer in 
$\ct'_\cc$ has order $i$. We have an isomorphism $\ov{Z(s')}^F@>\si>>\ct'_\cc$ (where 
$Z(s')$ is the centralizer of $s'$ in $PSO_Q$) induced by the homomorphism $g\m t$,
$Z(s')^F@>>>\ct'{}^F$ where $t=\dg s\dg\i s\i$; here $\dg\in G^{*F}$ is such that 
$\p(\dot g)=g$.) We can also view $\tz_{i,\cc}$ as the number of elements in
$\fU_F(Z(s))=\fU_F(Z(s')^0)$ whose stabilizer in the natural
$\ov{Z(s')}^F$ action on $\fU_F(Z(s')^0)$ has cardinal $i$. Hence, according to 
4.6(a) we have
$$|\ov{Z(s')}^F|\i\sum_{i\ge1}i^2\tz_{i,\cc}=
|\Irr_{\ov{Z(s')}^F}(\fU_F(Z(s')^0))|.$$
Using this and (c) we deduce
$$|\ct'{}^F|\i\sum_{i\ge1}i^2\tz_{i,\cc'}=
|\Irr_{\ov{Z(s')}^F}(fU_F(Z(s')^0))|.\tag d$$
We now sum the sum the identities (d) over all $\cc'\in\fs(PSO_Q^F)$. By (b), the sum of
the left hand sides becomes $|\ct'{}^F|\i\sum_{i\ge1}i^2\tz_i$. The sum of the right 
hand sides becomes $\ha_N^0$ (see 4.17(a)). Since $|\ct'{}^F|=|C^*_0|=|C^F|$, we 
obtain
$$|C^F|\i\sum_{i\ge1}i^2\tz_i=\ha_N^0$$
and using 5.3(a):
$$|C^F|\i\sum_{i\ge1}i^2\tz_i=|\Irr(\tSO_Q^F)|.\tag e$$
Applying \cite{\AST, 9(c)} to $B=G^F,A=\tSO_Q^F$ we obtain
$$|C^F|\i(\tz_{i,0}+4(\tz_{2,0}+y)+16(\tz_{4,0}-y))=|\Irr(\tSO_Q^F)|$$
where $y$ is the number of representations in $\Irr(G^F)$ whose restriction to
$\tSO_Q^F$ is not multiplicity free. Substracting this from (e) we obtain
$-12|C^F|\i y=0$ that is, $y=0$. In other words,

(f) {\it any irreducible representation of $G^F$ is multiplicity free when restricted to 
$\tSO_Q^F$.}
\nl
Now the proof of Theorem 5.2 for a general $G$ can be carried out as in \cite{\AST} using 
the multiplicity one statement (f).

\enddocument